\documentclass[11pt]{amsart} 

\usepackage{fullpage}

\usepackage{url}

\usepackage{amssymb}

\usepackage{cite}

\usepackage{graphicx}

\usepackage[usenames]{color}

\usepackage{accents}

\renewcommand{\a}{\alpha}

\newcommand{\g}{\gamma}
\renewcommand{\d}{\delta}
\newcommand{\D}{\Delta}
\newcommand{\e}{\varepsilon}
\newcommand{\f}{\varphi}
\newcommand{\s}{\sigma}
\renewcommand{\S}{\Sigma}

\renewcommand{\l}{\lambda}

\renewcommand{\O}{\Omega}

\newcommand{\cF}{{\mathcal F}}
\newcommand{\cC}{{\mathcal C}}
\newcommand{\cM}{{\mathcal M}}

\newcommand{\cS}{{\mathcal S}}
\newcommand{\cB}{{\mathcal B}}

\newcommand{\cE}{{\mathcal E}}
\newcommand{\cU}{{\mathcal U}}

\newcommand{\cP}{{\mathcal P}}

\newcommand{\cH}{{\mathcal H}}
\newcommand{\cD}{{\mathcal D}}
\newcommand{\cR}{{\mathcal R}}

\newcommand{\cI}{\mathcal I}

\newcommand{\bR}{\mathbb R}

\newcommand{\be}{\begin{equation}}
\newcommand{\ee}{\end{equation}}

\newcommand{\tr}{\mathrm{tr}}

\newcommand{\beaa}{\begin{eqnarray*}}
\newcommand{\bea}{\begin{eqnarray}}
\newcommand{\beal}[1]{\begin{eqnarray}\label{#1}}
\newcommand{\bean}{\begin{eqnarray}\nonumber}
\newcommand{\beadl}[1]{\begin{deqarr}\label{#1}}
\newcommand{\eeadl}[1]{\arrlabel{#1}\end{deqarr}}
\newcommand{\eeal}[1]{\label{#1}\end{eqnarray}}
\newcommand{\eead}[1]{\end{deqarr}}
\newcommand{\eea}{\end{eqnarray}}
\newcommand{\eeaa}{\end{eqnarray*}}

\newcommand{\Ric}{\operatorname{Ric}}

\renewcommand{\to}{\rightarrow}

\renewcommand{\exp}{\operatorname{exp}}

\DeclareMathOperator{\area}{area}

\DeclareMathOperator{\Image}{Im}

\DeclareMathOperator{\vol}{vol}
\DeclareMathOperator{\interior}{int}
\DeclareMathOperator{\dist}{dist}
\DeclareMathOperator{\Ker}{Ker}
\DeclareMathOperator{\Coker}{Coker}
\DeclareMathOperator{\Imm}{Imm}
\DeclareMathOperator{\Emb}{Emb}

\renewcommand{\phi}{\varphi}
\renewcommand{\epsilon}{\varepsilon}
\renewcommand{\hat}{\widehat}

\renewcommand{\>}{\rangle}

\newcommand{\dm}{{\partial M}}

\newcommand{\w}{\widetilde}

\theoremstyle{plain}
\newtheorem{theorem}{Theorem}[section]

\newtheorem{remark}[theorem]{Remark}

\newtheorem{lemma}[theorem]{Lemma}
\newtheorem{sublemma}[theorem]{Sub-Lemma}

\newtheorem{proposition}[theorem]{Proposition}
\newtheorem{corollary}[theorem]{Corollary}

\newtheorem{conjecture}[theorem]{Conjecture}

\theoremstyle{definition}

\def\endproof{\qed \medskip}
\def\blacksquare{\hbox to .60em {\vrule width .60em height .60em}}

\numberwithin{equation}{section}

\begin{document}

\title[ ]{Embeddings, immersions and the Bartnik quasi-local mass conjectures}

\author{Michael T. Anderson}
\address{Dept. of Mathematics, 
Stony Brook University,
Stony Brook, NY 11790}
\email{anderson@math.sunysb.edu}

\author{Jeffrey L. Jauregui}
\address{Dept. of Mathematics,
Union College,
Schenectady, NY 12308}
\email{jaureguj@union.edu}

\begin{abstract}
Given a Riemannian 3-ball $(\bar B, g)$ of non-negative scalar curvature, Bartnik conjectured that $(\bar B, g)$ admits an asymptotically flat (AF) 
extension (without horizons) of the least possible ADM mass, and that such a mass-minimizer is an AF solution to the static vacuum 
Einstein equations, uniquely determined by natural geometric conditions on the boundary data of $(\bar B, g)$. 

  We prove the validity of the second statement, i.e.~such mass-minimizers, if they exist, are indeed AF solutions of the static vacuum 
equations. On the other hand, we prove that the first statement is not true in general; there is a rather large class of bodies $(\bar B, g)$ for 
which a minimal mass extension does not exist.
\end{abstract}

\maketitle

\setcounter{section}{0}
\setcounter{equation}{0}

\section{Introduction}
\setcounter{equation}{0}

A fundamental problem in general relativity is the formulation of a ``suitable'' definition of quasi-local mass (cf. \cite[Problem 1]{Penrose}). 
To motivate this concept, consider for instance a time-symmetric, asymptotically flat (AF) initial data set $(M,g)$ without boundary for the Einstein equations, 
i.e.~a Riemannian 3-manifold viewed as a totally geodesic spacelike hypersurface in a Lorentzian (3+1)-dimensional spacetime. Assuming the spacetime 
obeys the dominant energy condition, the submanifold $(M,g)$ has non-negative scalar curvature. The quasi-local mass of a compact 
region $\Omega \subset (M,g)$ should be a real number that represents the mass contained within $\Omega$.

Many definitions of quasi-local mass have been put forth in the last several decades, though we make no attempt here to give a comprehensive 
history, see \cite{Sz} for an excellent review. Some of the ``classical'' examples include the Hawking mass \cite{H}, the Brown--York mass \cite{BY}, and 
the Bartnik mass \cite{Ba1}. More recently, Wang--Yau proposed a very interesting definition that generalizes the approach of 
Brown--York \cite{WY}.

In this paper we are interested in the Bartnik mass, whose setup we now recall. Let $\O$ be a smooth 3-manifold, with boundary, 
diffeomorphic to the closed 3-ball $\bar B$ in $\bR^3$, and let $g_{\O}$ be a Riemannian metric on $\O$ with non-negative scalar 
curvature. The Bartnik mass was originally defined as
\be \label{bm}
m_{B}(\O, g_{\O}) = \inf_{g} \{m_{ADM}(g)\},
\ee
where the infimum is taken over the set of smooth AF metrics $g$ on $\bR^3$ such that $(\Omega, g_{\O})$ embeds 
isometrically into $(\bR^3, g)$, and $(\bR^3, g)$ has non-negative scalar curvature and contains no horizons \cite{Ba1}. (Note that a smooth AF 3-manifold of non-negative scalar curvature with no horizons is diffeomorphic to $\bR^3$ \cite{MSY}.) Bartnik defined a 
horizon to be a stable minimal 2-sphere, but a number of variants have since been considered in the literature. Among these, we will take a 
horizon to be an immersed compact minimal surface that surrounds $\Omega$; this choice is discussed further in Section 2. 

The Bartnik mass satisfies many of the generally desired properties of a quasi-local mass (cf. \cite{Ba1}). For instance, $m_{B}(\O, g_{\O})$ 
is non-negative, by the positive mass theorem \cite{SY}, \cite{W}. Furthermore, if $(\O, g_{\O})$ is isometric to a smooth region in Euclidean space 
$(\bR^3, g_{Eucl})$, then $m_{B}(\O, g_{\O}) $ vanishes. Bartnik conjectured that  the converse holds (``strict positivity of $m_B$''), 
i.e., if $m_B(\O, g_{\O})=0$ then $(\O,g_{\O})$ is a Euclidean region. A key result of Huisken and Ilmanen \cite{HI} shows that 
if  $m_B(\O, g_{\O})=0$, then $(\O, g_{\O})$ is {\em locally flat}, i.e.~locally isometric to Euclidean space (although this result applies to a slightly different definition of Bartnik mass). The Bartnik mass also 
enjoys monotonicity (i.e.~a region contained in $(\O, g_{\O})$ cannot have a greater value of $m_B$; this follows from the definition), 
and the Bartnik mass limits to the ADM mass for an exhausting sequence of large balls in an AF manifold of nonnegative scalar curvature \cite{HI}. 
The most fundamental open questions regarding the Bartnik mass are to determine under which general conditions the infimum in \eqref{bm} is 
achieved, to understand the structure of the space of such minimizers and to describe the behavior of the corresponding 
mass functional on the space of minimizers. Before proceeding further, we recast the Bartnik mass in a slightly different 
manner, by focusing on the role played by the boundary geometry on the two sides of $\partial \Omega$.

\medskip

  For a pair $(\O, g_{\O})$ as above, let $\g = g_{\O}|_{T(\partial \O)}$ be the induced metric on $\partial \O \cong S^2$, and let $H$ be the 
mean curvature of $\partial \O$, (with respect to the unit outward normal, i.e. positive for round spheres in $\mathbb{R}^3$). The pair $(\g, H)$ will be called the (geometric) 
\emph{Bartnik boundary data} of $(\O, g_{\O})$. More generally, any pair $(\g,H)$, where $\g$ is a smooth Riemannian 
metric on $S^2$ and $H$ is a smooth function on $S^{2}$, will be called Bartnik boundary data.

  Bartnik pointed out that a minimizer $g$ of \eqref{bm} would only be expected to be  Lipschitz along the ``seam'' $\partial \Omega$, 
obeying the boundary conditions \cite{Ba2}, \cite{Ba3}
\be \label{bcont}
\g_{\partial \O} = \g_{\dm}, \ \ \ \ H_{\partial \O} = H_{\dm},
\ee
where $M$ is the closure of the complement of the embedded image of $\Omega$ in $\bR^3$ and $H_{\partial M}$ is the mean curvature of $\partial M$ with respect to the unit normal pointing into $M$. The significance of matching the mean curvatures on 
both sides is that it assures the scalar curvature is distributionally non-negative across the seam. The scalar curvature is also well-known 
to be distributionally non-negative if
\be \label{bcont_weak}
\g_{\partial \O} = \g_{\dm}, \ \ \ \ H_{\partial \O} \geq H_{\dm},
\ee
are satisfied; we discuss this point further in Remark \ref{remark_BCs}. The boundary condition \eqref{bcont_weak} was also 
considered by Miao \cite{Mi1} and Shi--Tam \cite{ST}. Thus, we consider the following reformulation of the Bartnik mass. Fix $M$ as a 
smooth manifold-with-boundary diffeomorphic to the closure of $\bR^3 \setminus \bar B$, and consider the space 
$\cP(M)$ of smooth, AF Riemannian metrics $g$ on $M$ (see \eqref{eqn_AF}) with non-negative scalar curvature, with $\accentset{\circ}\cP(M)$ being the subset 
such that $(M,g)$ contains no horizons (defined as above: immersed compact minimal surfaces that surround $\Omega$). We call $g \in \cP(M)$ an \emph{admissible extension} of a region 
$(\Omega, g_{\Omega})$ as above if \eqref{bcont_weak} holds. We take as our definition of the Bartnik mass:
\be \label{bm2}
m_{B}(\O, g_{\O}) = \inf \{m_{ADM}(g): g \in \accentset{\circ} \cP(M) \ {\rm is \ an \ admissible \ extension \ of} \ (\O, g_{\O})\}.
\ee 
One might also consider the mass defined by the equality condition \eqref{bcont}. Both of these versions have previously appeared in the literature. For further discussion on the numerous variations in the definition of Bartnik mass, and some progress on reconciling them, see \cite{J3}, \cite{Mc3}.

   These three definitions, based on \eqref{bm}--\eqref{bcont_weak}, all require a precise choice among the various possible definitions of horizon. 
A major reason a horizon is defined here to be a surrounding minimal surface (as opposed to an arbitrary minimal surface in $M$) is that 
$\accentset{\circ} \cP(M)$ is then open in $\cP(M)$, cf.~Lemma \ref{lemma_horizons} below, so that this condition is stable. (This is unknown for other 
definitions of the horizon condition.)

  Regarding then the boundary conditions \eqref{bm}--\eqref{bcont_weak} themselves, we prove in Theorem \ref{thm_BC} below that if a minimizer 
subject to \eqref{bcont_weak} exists in $\accentset{\circ} \cP(M)$, it necessarily satisfies \eqref{bcont} (cf.~also prior work of Miao on this issue \cite{M1}). 
This result strongly suggests the two definitions of Bartnik mass based on \eqref{bcont} and \eqref{bcont_weak} are equivalent and also 
very likely equivalent to \eqref{bm}, cf.~Remark \ref{remark_BCs} and \cite{Mc3}, \cite{J3}. Henceforth, we adopt \eqref{bm2} as the definition of the Bartnik mass.

\medskip 

  The following three conjectures are due to Bartnik; they are discussed in \cite{Ba1}, \cite{Ba2} and in most detail in \cite{Ba3}. 
  
\medskip

\noindent  {\bf Conjecture I.} Any region $(\O, g_{\O})$ with $\Omega \cong \bar B$, $H_{\partial \O}>0$, and $g_{\O}$ of non-negative scalar curvature admits an admissible extension in $\accentset{\circ} \cP(M)$.

\smallskip 

  Thus, conjecturally, any metric of non-negative scalar curvature on a ball, with positive boundary mean curvature, can be extended to an AF manifold with non-negative 
scalar curvature, where the extension has no horizons and \eqref{bcont_weak} is satisfied. (The hypothesis of positive boundary mean 
curvature is imposed because if, for instance, $H_{\partial \O}$ were negative everywhere, then any AF extension would contain a 
horizon.) This general extension conjecture essentially appears in \cite[Problem 1]{Ba3}. It implies that any region ($\O$, $g_{\O}$) as 
above has a well-defined Bartnik mass \eqref{bm2}.

  Conjecture I is known as the Bartnik extension conjecture and remains open in general (even allowing extensions in $\cP(M)$). 
Further discussion of the conjecture and some partial results are given in Section 3. 

\medskip 

\noindent  {\bf Conjecture II.} For any region $(\O, g_{\O})$ with $\Omega \cong \bar B$, $H_{\partial \O} > 0$, and $g_{\O}$ of non-negative scalar curvature, there exists an admissible 
extension $g \in \accentset{\circ} \cP(M)$ realizing the Bartnik mass \eqref{bm2}. Moreover, $g$ satisfies the boundary conditions \eqref{bcont}.

\smallskip 

Conjecture II is known as the Bartnik mass-minimization conjecture. Bartnik \cite{Ba3}, \cite{Ba4} developed a heuristic program 
suggesting that a metric $g$ realizing the Bartnik mass \eqref{bm2} on $M$ yields an asymptotically flat (AF) solution of the static vacuum 
Einstein equations, i.e.~there is a potential function $u: M \to \bR$, with $u \to 1$ at infinity, such that 
\be \label{stat}
u\Ric_{g} = D^{2}u, \ \ \D u = 0.
\ee
(Moreover, the AF potential function $u$ is expected to be positive.) 
This has been partially verified, using quite different methods, by Corvino \cite{C1}, \cite{C2}, cf.~Remark \ref{corv} for further 
discussion. We give a full proof of this proposal, thus completely implementing Bartnik's program:

\begin{theorem} \label{thm1.1} 
For a region $(\Omega, g_\Omega)$ as above with $H_{\partial \O}>0$, a metric $g$ minimizing the Bartnik mass \eqref{bm2} admits an AF potential function $u>0$ such that $(g, u)$ is an AF solution 
of the static vacuum Einstein equations \eqref{stat}. Moreover, such a minimizer $g$ satisfies \eqref{bcont}. 
\end{theorem}

   We refer to Theorem \ref{thm_mass_min_static} and Theorem \ref{thm_BC} for further details. 

\medskip 

\noindent   {\bf Conjecture III.} For any geometric Bartnik boundary data $(\g, H)$ on $S^{2}$, with $H > 0$, there exists a unique extension 
$g \in \accentset{\circ} \cP(M)$ of $(\g, H)$ and a function $u>0$ with $u \to 1$ at infinity, such that the pair $(g,u)$ is an AF solution of the static vacuum Einstein equations \eqref{stat}.

\smallskip 

 Conjecture III is known as the Bartnik static metric extension conjecture. 

\medskip 

  In addition to the horizon issue, the assumption $H > 0$ in Conjectures II and III is made due to the black hole uniqueness 
theorem, cf.~\cite{I}, \cite{BM}, and also \cite{Mi2}. Namely, the data $(\g, 0)$ are boundary data of a static vacuum solution 
only for $\g$ a round, constant curvature metric on $S^{2}$, realized by the family of Schwarzschild metrics. Thus Conjectures 
II and III are well-known to fail for $H = 0$ boundary data. 

  It is clear that Conjectures II and III each imply Conjecture I.  Using Theorem \ref{thm1.1}, Conjecture II implies the existence 
part of Conjecture III for the special case of boundary data $(\g, H)$ obtained from a region 
$\Omega$ with non-negative scalar curvature and $H > 0$. On the other hand, even for this special class of boundary data, Conjecture III 
does not imply Conjecture II, since all mass-minimizing sequences for a given body $\O$ may fail to converge to a limit. 
As discussed in Proposition \ref{prop_critical}, the static vacuum solutions given in Conjecture III are precisely the critical points of the ADM mass 
$m_{ADM}$ (with fixed boundary conditions), but it is not clear that these are minimizers. If Conjecture II holds, so minimizers exist, 
then the uniqueness of Conjecture III would imply that all critical points are minimizers. 

\medskip

   Given this background, the main purpose of this work is to prove that Conjecture II is not true in the generality 
stated, so that further hypotheses are required to maintain its validity (see the discussion at the end of Section 5). As discussed below, 
similar remarks apply to Conjecture III. This failure is related to the degeneration of the exterior manifold-with-boundary structure on $M$, 
given control on the boundary data in \eqref{bcont} or \eqref{bcont_weak}. This is most simply described in the passage 
from embedded spheres to immersed spheres in $\bR^3$, which we now discuss.

  Let $\Imm(\bar B, \bR^{3})$ be the space of smooth immersions 
$$F: \bar B \to \bR^{3},$$
of the closed 3-ball $\bar B \subset \bR^{3}$; thus $F$ extends to an immersion of an open neighborhood of $\bar B$. 
Let $\cF \subset \Imm(\bar B, \bR^{3})$ denote the subspace of immersions that restrict to embeddings on the interior of $\bar B$ and on 
which the self-intersection set $Z$ of $F|_{\partial \bar B} = F|_{S^{2}}$ consists of a finite, nonzero number of double points. Thus, there is a 
finite set $Z \subset S^2$ given as a disjoint union $(\cup z_{i}) \cup (\cup z_{i}') $ such that $F$ is injective on $\bar B \setminus Z$,  $F(z_i) = F(z_i')$ 
for each $i$, and $F(z_i) \neq F(z_j)$ for $i \neq j$. For $F \in \cF$, the set $F(\bar B)$ is not a smooth region in $\bR^{3}$. However, the pullback 
$(\bar B, F^*(g_{Eucl}))$ is obviously a smooth, locally flat Riemannian manifold with boundary. It is easy to see that $\cF$ provides a large, 
infinite-dimensional space of such locally flat domains. We also remark that there is a large class of immersions $F \in \cF$ such that 
$(\bar B, F^*(g_{Eucl}))$ has positive boundary mean curvature (consider, for example, a torus in $\mathbb{R}^3$ with positive mean curvature, appropriately cut with spherical caps glued back in, tangent at their respective poles --- see \cite[Figure 1]{AK}).

\begin{theorem}  \label{thm_conj_II} 
Conjecture II is false for any region $(\bar B, F^{*}(g_{Eucl}))$ for $F \in \cF$ as above. In particular, there is no admissible extension 
of $(\bar B, F^{*}(g_{Eucl}))$ whose ADM mass attains the Bartnik mass (which equals zero).
\end{theorem}

  In particular, this also shows that strict positivity of the Bartnik mass fails, i.e.~ the result of Huisken--Ilmanen \cite{HI} that 
$m_B(\O, g_{\O})=0$ implies local flatness is optimal. This is because the proof of Theorem \ref{thm_conj_II} will 
show that $(\bar B, F^{*}(g_{Eucl}))$ has zero Bartnik mass and does not \emph{embed} isometrically in Euclidean 3-space, 
cf.~also Remark \ref{rem_strict_pos}. 

   Note that for $F \in \cF$, there is a sequence of embeddings $F_{i}$ of the closed 3-ball into $\bR^3$ with $F_{i} \to F$ smoothly, 
with the corresponding embedded spheres $F_{i}(S^{2})$ converging smoothly to an immersed sphere. In particular the 
class $\cF$ of immersions is at the boundary of the space of embeddings. Of course Conjecture II holds for regions 
$\O$ isometrically embedded in $\bR^{3}$. 

  As noted above, the pulled-back Euclidean metrics $F_{i}^{*}(g_{Eucl})$ converge smoothly to a limiting smooth flat metric 
on the abstract $3$-ball $\bar B$ with limit boundary data $(\g, H)$. However, the flat metrics on the complementary 
manifolds $M_{i} = \bR^{3}\setminus F_i(\bar B)$ degenerate in the limit. It is a priori possible that there is a distinct 
sequence of (non-flat) admissible extensions $g_{i}$ of the boundary data $(\g, H)$ with $m_{ADM}(g_{i})$ converging 
to the infimum of the mass of such extensions, which do not degenerate and so give a limit realizing the Bartnik mass. 
The main content of Theorem \ref{thm_conj_II} is to prove that in fact this does not occur. 

  We conjecture that this phenomenon is quite general, i.e. that Conjecture II is false for any domain $(\bar B, F^{*}(g_{Eucl}))$ obtained 
from an immersion $F: \bar B \to \bR^3$ that is not an embedding (even if $F$ is not at the boundary of the space of 
embeddings), cf.~Conjecture \ref{conj_false}.
  
\medskip   
  
  A version of the discussion above also holds with respect to Conjecture III. Namely, let $\cE^{m,\a}$ be the 
moduli space of AF static vacuum solutions $(g, u)$, $u > 0$, on $M = \bR^{3}\setminus B$. The moduli space $\cE^{m,\a}$ 
is the space of all static vacuum metrics $(g, u)$ which are $C^{m,\a}$ smooth up to $\dm$, modulo the action of the 
$C^{m+1,\a}$ diffeomorphisms ${\rm Diff}_{1}^{m+1,\a}(M)$ of $M$ equal to the identity on $\dm$ (and asymptotic to the 
identity at infinity). It is proved in \cite{A1} (cf.~also \cite{AK}) that $\cE^{m,\a}$ is a smooth Banach manifold, and moreover 
the map to Bartnik boundary data 
\be \label{PiBa}
\Pi_{B}: \cE^{m, \a} \to Met^{m,\a}(S^{2})\times C^{m-1,\a}(S^{2}),
\ee
$$\Pi_{B}(g,u) = (\g, H),$$
is a smooth Fredholm map, of Fredholm index 0. Here, $Met^{m,\a}(S^{2})$ is the space of $C^{m,\a}$ Riemannian metrics on $S^2$
with the $C^{m,\a}$ topology. 

  Now consider the map $\Pi_{B}^+$, the  restriction of $\Pi_B$ to the open subspace $\cE_{+}^{m,\a}$ of static vacuum metrics 
with $H > 0$ at $\dm$:  
\be \label{PiB+}
\Pi_{B}^+: \cE_{+}^{m,\a} \to Met^{m,\a}(S^{2})\times C_{+}^{m-1,\a}(M).
\ee
Conjecture III is equivalent to the statement that $\Pi_{B}^+$ in \eqref{PiB+} is a bijection. However, it is proved in 
\cite{AK} that $\Pi_{B}^+$ is not a homeomorphism; in fact the inverse map to $\Pi_{B}^+$, if it exists, is not continuous 
in general. The failure of the homeomorphism property is closely related to the behavior of $\Pi_{B}^+$ at the boundary of the 
space of (flat) embeddings within the larger space of immersions, discussed above in connection with Conjecture II. 

\medskip 

  Theorem \ref{thm_conj_II} shows that a major obstacle in establishing the validity of Conjecture II is controlling the behavior of 
mass-minimizing sequences arbitrarily close to the boundary $\dm = \partial \O$, given control on the Bartnik boundary 
data $(\g, H)$, so that the manifold-with-boundary structure of $M$ does not degenerate. A similar difficulty arises 
in proving Conjecture III; for example, it is much simpler to control the behavior of sequences of static vacuum solutions in the 
interior of $M$ (away from $\dm$) compared with controlling the behavior near the boundary; see for example the analysis in \cite{A0}. 
We expect a similar phenomenon for more general mass-minimizing sequences. 

\medskip 

  In contrast to the negative results above on Conjectures II and III, we present positive evidence for the validity 
of Conjecture I in Section 3. We prove in Proposition \ref{prop_H} that if the boundary data $(\g, H)$ admit an extension to an AF metric 
of non-negative scalar curvature, then so do $(\g, \w H)$, for any $\w H \geq H$. Combining this with previous results in \cite{MS} and 
\cite{A2} leads to the verification of Conjecture I for a wide variety of boundary data $(\g, H)$, although without addressing the issue of 
horizons.

\medskip 

  The contents of the paper are briefly as follows. In Section 2, we discuss the various possible definitions of horizon 
as well as the boundary conditions \eqref{bcont}--\eqref{bcont_weak}, and the relations of minimizers of the Bartnik mass 
with the static vacuum Einstein equations. The main results are Theorems \ref{thm_mass_min_static} and \ref{thm_BC}, which imply Theorem \ref{thm1.1}. In Section 3 we discuss Conjecture I and present new evidence for its validity in general. Section 4 is devoted to the proof of 
Theorem \ref{thm_conj_II}, while Conjecture III is discussed further in Section 5. We note that although the topics of these sections are 
of course inherently related, the sections themselves are essentially independent of each other. 

\medskip 

\noindent \textit{Acknowledgments:}
   This work began at the conference ``Static metrics and Bartnik's quasi-local mass conjecture" in May 2016 at the 
Universit\"at T\"ubingen. We are grateful to the University for financial support and in particular to Carla Cederbaum, 
for providing the opportunity to initiate this collaboration. We also extend our thanks to Zhongshan An and the referees for 
their careful reading of the paper and insightful questions and comments. M.A.~was partially supported by NSF Grant DMS 1607479.

\section{Mass minimizers and the static vacuum Einstein equations}
\setcounter{equation}{0}

  In this section, we discuss relations between the various notions of Bartnik mass from the Introduction and their relations 
with the static vacuum Einstein equations. 

\medskip 

  Starting with an idea suggested by Brill--Deser--Fadeev in \cite{BDF}, Bartnik in \cite{Ba3},\cite{Ba4} presented a heuristic 
argument that critical points of the mass on the space of solutions of the (time-symmetric) 4-dimensional vacuum 
Einstein constraint equations, with fixed boundary data $(\g, H)$, should be  given by solutions of the static vacuum Einstein 
equations. This strongly suggested that minimizers of the Bartnik mass (with respect to a suitable horizon condition) should 
then also be static vacuum Einstein solutions. Some recent work along these lines has also been carried out by 
McCormick \cite{Mc1}, \cite{Mc2}.

  The main results of this section are a full proof of Bartnik's proposal for Bartnik mass minimizers, cf.~Theorems \ref{thm_mass_min_static} and \ref{thm_BC}. In addition, 
Theorem \ref{thm_BC} shows that a Bartnik mass minimizer defined according to \eqref{bcont_weak} actually satisfies \eqref{bcont}, 
leading to a corresponding strong monotonicity result in Corollary \ref{mincor}. We refer the reader to Remark \ref{corv} for relevant prior results.

\medskip 

  Throughout, $M$ will be a smooth 3-manifold with boundary, diffeomorphic to $\bR^{3}\setminus B$, where $B$ is an open ball. 
A $C^{m,\a}$ Riemannian metric $g$ on $M$ (i.e.~$C^{m,\a}$ up to and on $\dm$) will be called asymptotically flat (AF) if
\be \label{eqn_AF}
g_{ij} - \d_{ij} \in C_{\d}^{m,\a}(M),
\ee
where $C_{\d}^{m,\a}(M)$ is the weighted H\"older space of functions on $M$ that decay to 0 at a rate $r^{-\d}$ with 
$k^{\rm th}$ derivatives decaying at the rate $r^{-\d - k}$, $k \leq m$, and with appropriately weighted H\"older $\alpha$-seminorms of the $m$th derivatives bounded, cf.~\cite{CD} for example. The rate 
$\d$ is fixed throughout and assumed to satisfy
$${\tfrac{1}{2}} < \d < 1.$$
We also fix any $m \geq 3$ and $\a \in (0,1)$. 

  Given a fixed $\d$ as above, let $\cP^{m,\a}(M) = \cP_{\d}^{m,\a}(M)$ be the space of AF metrics $g$ on $M$ with non-negative 
scalar curvature $s = s_g$. Recall that the ADM mass $m_{ADM}$ of $g \in \cP^{m,\a}(M)$ is only defined \cite{Ba0} for metrics with 
\be \label{l1}
s \in L^{1}(M).
\ee
The Bartnik mass \eqref{bm2} is then obtained by minimizing the ADM mass on $\cP^{m,\a}(M)$ subject to the boundary 
conditions \eqref{bcont_weak} on $\dm$ and subject to the no-horizon condition. Alternatively, one might consider 
minimizing the mass subject to the stronger condition \eqref{bcont}. 

  As mentioned in the Introduction, there are several notions of horizon appearing in the literature 
without a current general consensus. The most strict condition is that $(M, g)$ has no immersed 
compact minimal surfaces; let $m_{B}^{s}$ denote the corresponding Bartnik mass. Variations of this condition such as 
no stable compact minimal surfaces or no embedded compact minimal surfaces have also been considered. In some cases, 
the minimal surfaces are required to be topological spheres.

  A somewhat weaker condition (and the one that we adopt) is that there are no immersed compact minimal surfaces surrounding $\dm$ in $M$, 
i.e.~ any path from $\dm$ to infinity must pass through the surface. (Again one might consider variations such as no 
stable or no embedded surrounding compact minimal surfaces).  Let $m_{B}^{w}$ denote the corresponding Bartnik mass; then one clearly has 
\be \label{ws}
m_{B}^{w} \leq m_{B}^{s}.
\ee
The same relation holds with respect to the weaker and stronger boundary conditions \eqref{bcont_weak} 
and \eqref{bcont}, respectively. 

  Moreover, a third definition was suggested by Bray \cite{Br}, requiring that $\dm$ be outer-minimizing in $(M, g)$. 
This version of the mass will be discussed briefly in Section 5, but not used before then. 

   One of the main reasons for preferring the weaker condition is the following stability result. Let 
\be \label{P0}
\accentset{\circ} \cP ^{m,\a}(M) \subset \cP^{m,\a}(M)
\ee
be the subset of metrics that have no immersed minimal surface surrounding $\dm$.

\begin{lemma}
\label{lemma_horizons}
$\accentset{\circ} \cP^{m,\a}(M)$ is an open subset of $\cP^{m,\a}(M)$.
\end{lemma}

\noindent {\bf Proof:} We show that the complement is closed. Let $\{g_i\}$ be a sequence in $\cP^{m,\a}(M)\setminus \accentset{\circ} \cP^{m,\a}(M)$ 
converging to some $g \in \cP^{m,\a}(M)$. Each $g_i$ is an AF metric on $M$ such that $(M,g_i)$ contains an immersed minimal surface 
$\S_i$ surrounding $\dm$. 
The unbounded component $\hat M_i$ of $M \setminus \S_i$ is then AF; one may then minimize the area functional for surfaces in $\hat M_i$ homologous to infinity. Since $\partial \hat M_i$ together with a large sphere $S$ near infinity (independent of $i$), serve as well-defined barriers, it follows from well-known 
results of Meeks--Simon--Yau \cite{MSY} that $\hat M_i$ contains a minimal surface $\hat \S_i$ 
that has the least area among surfaces surrounding $\partial \hat M_i$. 

In particular, $\hat \S_i$ is stable. Further, the area of $\hat \S_i$ 
with respect to $g_i$ is uniformly bounded, since $\hat \S_i$ has less $g_i$-area than $S$, and $\area_{g_{i}}(S) \to \area_{g}(S)$. 
Using the well-known curvature estimates of Schoen, it is then standard, (cf.~\cite{CM} for example) that a subsequence of $\hat \S_i$ 
converges to a stable minimal surface $\S$ in $(M,g)$. Clearly $\S$ surrounds $\partial M$, so that $g \in \cP^{m,\a}(M)\setminus 
\accentset{\circ} \cP^{m,\a}(M)$.

{\endproof}

   It is mainly for this stability behavior that we choose a horizon to be a  surrounding minimal surface. (Such stability is unknown for 
other definitions of horizon in an AF manifold with boundary.) A further reason is that static vacuum Einstein metrics have no horizons 
in this sense (except for Schwarzschild metrics), by a result of Miao \cite{Mi2}. (This is again unknown for general minimal surfaces). We prove Theorem \ref{thm_conj_II} for the version of the Bartnik mass we adopt (i.e., for horizons as surrounding minimal surfaces and for the weaker boundary condition 
\eqref{bcont_weak}). In Remark \ref{rmk_bcs}, we describe a possible modification to the proof that may allow for the stronger boundary condition \eqref{bcont}. Theorem \ref{thm_BC} below 
strongly suggests that with respect to $\accentset{\circ} \cP^{m,\a}(M)$ as in \eqref{P0}, the boundary conditions \eqref{bcont_weak} and 
\eqref{bcont} give equal Bartnik masses; this is less clear for the stronger definition $m_{B}^{s}$. 

  To summarize using the current notation, as in \eqref{bm2}, we set 
\be \label{bm21}
m_{B}(\O, g_{\O}) = m_{B}(\g, H) = \inf \{m_{ADM}(g): g \in \accentset{\circ} \cP^{m,\a}(M), g|_{\dm} = \g, H_{\dm} \leq H\}.
\ee   
Note that an immediate consequence of the definition is the following (weak) inverse monotonicity property: if $H' \leq H$, 
then 
\be \label{mono}
m_{B}(\g, H) \leq m_{B}(\g, H').
\ee
A strong monotonicity will be proved in Corollary \ref{mincor} below. 

\medskip 

  Returning to the discussion prior to \eqref{l1}, let $\cS^{m,\a}(M) = \cS_{\d}^{m,a}(M)$ be the space of pairs $(g, u)$, with $g$ 
a $C^{m,\a}_{\d}$ AF metric on $M$ and $u$ an AF function, i.e. $u - 1 \in C_{\d}^{m,\a}(M)$, so that $u \to 1$ at infinity. We write 
$\cS^{m,\a}(M) =Met_{AF}^{m,\a}(M)\times C_{AF}^{m,\a}(M)$. Data $(g,u)$ for which $u>0$ correspond to AF Lorentzian  metrics on 
$\cM = \bR \times M$ of the form
\be \label{N}
g_{\cM} = -u^{2}dt^{2} + g.
\ee
Metrics of 
the form \eqref{N} (or such pairs $(g,u)$) will be referred to as \emph{static}; this is not to be confused with other notions of 
static (e.g. static vacuum or Corvino's definition of static in \cite{C1}).

  Clearly $\cS^{m,\a}(M) = Met_{AF}^{m,\a}(M)\times C_{AF}^{m,\a}(M)$ is a smooth Banach manifold.   
Let $\cS_{+}^{m,\a}(M) \subset \cS^{m,\a}(M)$ be the subset such that 
$$s_{g} \geq 0.$$
Thus, $\cS_{+}^{m,\a}(M) = \cP^{m,\a}(M)\times C_{AF}^{m,\a}(M)$. Note that if
$(g, u) \in \partial (\cS_{+}^{m,\a}(M))$, then the scalar curvature of $g$ is nonnegative and vanishes at some point in $M$.  We point out that the condition \eqref{l1} is not assumed a priori on $\cS_{+}^{m,\a}$.

  Consider the Regge--Teitelboim Hamiltonian \cite{RT} in this setting: 
$$\cH: \cS^{m,\a}(M) \to \bR,$$
\be \label{cL}
\cH(g, u) = \int_{M}us dv_{g} - 16\pi m_{ADM}(g),
\ee
where $s = s_{g}$ is the scalar curvature of $g$.  If $u>0$, note that since $s_{g_{\cM}} = s_{g} - 2\frac{\D u}{u}$ and $dv_{g_{\cM}} = udv_{g}$, 
the first term gives the Einstein--Hilbert action on the 4-manifold $\cM$ modulo a divergence term (namely $-2\D u$). The reason 
for this modification of the Einstein--Hilbert action is to obtain a well-defined variational problem for the ADM Hamiltonian; 
we refer to \cite{RT} for details. 

  If $s \notin L^{1}(M)$, then the individual terms in \eqref{cL} are ill-defined although the combination is well-defined. 
Explicitly, following \cite{Ba4}, \eqref{cL} may be rewritten in the form 
\be \label{cL2}
\cH(g, u) = \int_{M}(u - 1)s dv_{g} - \int_{M}(\cR_{0} - s)dv_{g},
\ee
where $u - 1 \in C_{\d}^{m,\a}(M)$ and $\cR_{0}$ is the bulk integral for the mass, given by 
$$\cR_{0} = \left(\d_{0} \d_{0} g -  \D_{g_{0}}(\tr_{g_{0}} g)\right)\frac{dv_{g_{0}}}{dv_{g}},$$ 
where $g_{0}$ is any background metric agreeing with $g$ near $\dm$ and is Euclidean outside a compact set, and $\d_0$ 
is the corresponding divergence. By the divergence theorem, $\int_{M}\cR_{0}dv_{g} = 16\pi m_{ADM}(g)$, when the ADM mass 
is defined, cf.~\cite{Ba4}. Thus the Regge--Teitelboim Hamiltonian \eqref{cL2} is well-defined and a smooth functional on the 
full Banach manifold $\cS^{m,\a}(M)$.  Of course the definitions \eqref{cL} and \eqref{cL2} agree when $s \in L^{1}(M)$. 
 
  Let 
\be \label{S*}
S^{*}u = D^{2}u - (\D u) g - u\Ric,
\ee
be the formal $L^2$-adjoint of the linearization $s'=Ds_g$ of the scalar curvature. The static vacuum Einstein equations 
\eqref{stat} are equivalent to the following system for $(g, u) \in \cS^{m,\a}(M)$:
\be \label{stat1}
S^{*}u = 0, \ \ \D u = 0.
\ee
Note that static vacuum metrics are necessarily scalar-flat, $s = 0$, and are also necessarily in $\accentset{\circ} \cP^{m,\a}(M)$, 
i.e.~as noted above, have no horizons (except for Schwarzschild metrics), cf.~\cite{Mi2}. The relation 
$\D u = 0$ in \eqref{stat1} follows from $S^{*}u = 0$ by taking the trace, and using the Bianchi identity to establish constant scalar curvature; then asymptotic flatness gives $s=0$ and hence $\D u = 0$.

\medskip 

  The following result is essentially classical and is a version of results proved in \cite{RT}, \cite{FMM}, \cite{Ba4}; a simple proof 
in this notation is also given in \cite{AK}.  Let $N$ be the unit normal at $\dm$ pointing into $M$ and let $A$ be the $2^{\rm nd}$ 
fundamental form of $\dm$ in $M$. 

\begin{proposition} \label{l2grad}
The $L^{2}$-gradient of $\cH$ on $\cS^{m,\a}(M)$ is given by
\be \label{grad}
\nabla \cH = (S^{*}u + {\tfrac{1}{2}}us g, s, uA - N(u)\g, 2u)
\ee
in the sense that, if $(h, u')$ is any variation of $(g, u)$ inducing the variation $(h^{T}, H'_{h})$ of boundary data $(\g, H)$, 
then
\be \label{var}
d\cH_{(g, u)}(h, u', h^{T}, H'_{h}) = \int_{M}[\<S^{*}u + {\tfrac{1}{2}}usg, h\> + su'] + \int_{\dm}
[\<uA - N(u)\g, h^{T}\> + 2uH'_{h}].
\ee
\end{proposition}
(The volume forms associated with the metrics on $M$ and $\dm$ are omitted, to simplify the notation). 

  Let $\cS_{(\g,H)}^{m,\a}(M)$ be the space of static metrics with fixed Bartnik boundary conditions; 
thus $\cS_{(\g, H)}^{m,\a}(M)$ consists of pairs $(g, u) \in \cS^{m,\a}(M)$ with the metric $g$ having fixed boundary data 
equal to $(\g, H)$ at $\dm$. It is straightforward (using the implicit function theorem) to show that $\cS_{(\g,H)}^{m,\a}(M)$ is a smooth, closed Banach submanifold of 
$\cS^{m,\a}(M)$, for all choices of $(\g, H) \in Met^{m,\a}(S^{2})\times C^{m-1,\a}(S^{2})$. Tangent vectors to $\cS_{(\g,H)}^{m,\a}(M)$ 
are variations $(h,u')$ of $(g, u)$ such that $(h^{T}, H'_{h}) = (0,0)$ at $\dm$, where $h^{T}$ is the restriction of $h$ to $T(\dm)$ 
and $H'_{h}$ is the variation of the mean curvature in the direction of $h$.

  Proposition \ref{l2grad} thus shows that critical points of the Hamiltonian $\cH$ on $\cS_{(\g,H)}^{m,\a}(M)$ 
are given exactly by static vacuum Einstein metrics realizing the given boundary data $(\g, H)$. 

\medskip 

  In contrast, we show next that there are no critical points of the mass 
$$m_{ADM}: \cD \subset \cS_{(\g,H)}^{m,\a}(M) \to \bR,$$ 
where $\cD$ is the domain on which $m_{ADM}$ is well-defined, i.e. the subset of $\cS_{(\g,H)}^{m,\a}(M)$ consisting of metrics $g$ with integrable scalar curvature. Given $H$, let $\cS_{(\g,H^{\leq})}^{m,\a}(M)$ be the space 
of static metrics with  boundary metric $\g$ and mean curvature $\leq H$ at $\dm$. 

\begin{lemma} 
\label{lemma_grad_m}
For any $(g,u) \in \cS_{(\g,H)}^{m,\a}(M)$ for which $m_{ADM}$ is defined, one has 
\be \label{grad2}
(Dm_{ADM})_g \neq 0.
\ee
(i.e., there exists an appropriate variation $h$ of $g$ with $(Dm_{ADM})_g(h) \neq 0$). If, in addition, $(g,u)  \in \cS_{+}^{m,\a}(M)$, then $(Dm_{ADM})_g$ is non-vanishing in the directions of 
$\cS_{+}^{m,\a}(M)\cap \cS_{(\g,H)}^{m,\a}(M)$ (again in the sense that there exists an appropriate variation). Furthermore, if $(g,u)  \in \cS_{+}^{m,\a}(M)$ and
\begin{equation*}
s_{g} \not\equiv 0,
\end{equation*}
then there is an infinitesimal deformation $(h,0)$ of $(g,u)$ in the direction of $\cS_{+}^{m,\a}(M)\cap 
\cS_{(\g,H^{\leq})}^{m,\a}(M)$ such that 
\be \label{ml0}
(Dm_{ADM})_{g}(h) < 0,
\ee
so that there are metrics $g' \in \cS_{+}^{m,\a}(M)\cap \cS_{(\g,H^{\leq})}^{m,\a}(M)$ with $m_{ADM}(g') < m_{ADM}(g)$. 
\end{lemma}

\noindent {\bf Proof:} We use a well-known conformal argument, cf. \cite{C1} for example. Suppose $g$ is an AF metric and 
$\w g = v^{4}g$ is a conformal deformation of $g$, with $v >0 $ in $C^{m,\a}_{AF}$, so that $\w g$ is AF. The scalar curvatures 
of $\w g$ and $g$ are related by 
$$v^{5}\w s = -8\D v + sv.$$
Suppose the ADM mass $m$ of $g$ is defined, and that $\D v \in L^1(M)$. Then the ADM mass $\w m$ of $\w g$ is also defined, 
and a well-known formula (cf. \cite{Mi1}, eqn.~(46)) relating $m$ and $\w m$ reads 
\be \label{confmass}
\w m = m - \frac{1}{2\pi}\lim_{r \to \infty}\int_{S(r)}N(v)dA, 
\ee
where $N$ is the outward unit normal at the coordinate sphere $S(r)$, and $dA$ is the induced area form on $S(r)$ with respect to $g$. 

Let $\f$ be a superharmonic function on $(M, g)$, $\D \f \leq 0$, with $\f = 0$ in a neighborhood of $\dm$ and with $\f$ 
harmonic outside a large compact set, tending to a constant $-c < 0$ at infinity. It is easy to see such functions exist. Let $F_{t}: M \to M$ be a smooth family of diffeomorphisms equal to the identity near $\dm$ and equal to the map $x \to (1-tc)^{-2}x$ 
near infinity with $F_{0} = Id$. We apply \eqref{confmass} to the curve of metrics $g_{t} = F_{t}^{*}((1 + t\f)^{4}g)$.  (The diffeomorphisms are needed to put the curve $g_{t}$ in the space $Met_{AF}^{m,\a}(M)$). Note that $s_{g} \geq 0$ implies $s_{g_{t}} \geq 0$. 

   Taking the derivative of \eqref{confmass} and using the divergence theorem gives, for $r$ sufficiently large, 
$$m'_{h} =(Dm_{ADM})_g(h)=  -\frac{1}{2\pi}\int_{S(r)}N(\f) dA > 0,$$
for the variation $h = \partial_{t}g_{t}|_{t=0}$. (Note the diffeomorphisms $F_t$ may be neglected in this calculation, 
since the ADM mass is diffeomorphism invariant). Since $h$ preserves the boundary conditions, this proves the first two 
statements. 

  To prove the last statement, let $v \in C^{m,\a}_{AF}$ be the unique solution to the equation $-8\D v + s\rho v = 0$ with $v = 1$ 
on $\dm$ and $v \to 1$ at infinity, where $s\geq0$, $s \not\equiv0$, and $0 \leq \rho \leq 1$ is smooth and compactly supported, such that $\rho s \not\equiv 0$. In particular, $v \in C^{m,\a}_{AF}$. By the minimum principle, $v > 0$ on $M$, so that $\w g = v^{4}g$ 
is well-defined asymptotically flat, and $s_{\w g} \geq 0$. By the maximum principle, $v < 1$ in the interior of $M$ (since $\rho s$ is not identically zero) and $N(v) < 0$ at $\dm$. There is an $\e > 0$ such that the level set $v^{-1}(1 - \e)$ has a 
compact, regular component $L$ that is a closed surface in the AF end of $M$, enclosing the support of $\rho$. Thus $v$ is $g$-harmonic outside of $L$, equalling $1-\epsilon$ on $L$ and approaching 1 at infinity. By the Hopf maximum principle, $N(v)>0$ along $L$. Applying the divergence theorem in the region 
outside of $L$, it follows that 
$$\w m = m - \frac{1}{2\pi}\lim_{r\to\infty}\int_{S(r)} N(v)dA < m.$$
Moreover, one has $\w H = H + 4N(v) < H$ on $\partial M$, by the Hopf maximum principle. 

   One may also linearize this argument by choosing $v = v_{t}$ as above solving $-8\D v + ts\rho v = 0$, so that $-8\D v_{t} + sv_{t} = 
sv_{t} - tsv_{t} = sv_{t}(1-t) \geq 0$, i.e.~the conformally deformed metric has non-negative scalar curvature.  
Taking the derivative at $t = 0$ gives the result. 

{\endproof}

  Lemma \ref{lemma_grad_m} shows the special role accorded to the scalar-flat metrics on $M$. In the context of the 4-metric \eqref{N} on 
$\cM$, the equation 
$$s = 0,$$
on $(M, g)$ is exactly the set of vacuum Einstein constraint equations (since the $2^{\rm nd}$ fundamental form $K$ 
of $M$ in $\cM$ vanishes). Let  
$$\cC^{m,\a} = \{(g, u): s_{g} = 0\} \subset \cS^{m,\a}(M).$$
be the space of solutions of the constraint equations on $\cS^{m,\a}(M)$. Note that there is no condition on the lapse $u$ 
(except $u-1 \in C^{m,\a}_{\d}$). Of course $\cC^{m,\a}$ is a (small) subset of the full boundary $\partial(\cS_{+}^{m,\a}(M))$. 
Note also that since static vacuum metrics are scalar-flat, any critical point $(g, u)$ of $\cH$ on $\cS_{(\g,H)}^{m,\a}(M)$ must lie 
in $\cC^{m,\a}$. 

\medskip

To proceed further, we need to examine the smoothness of the spaces $\cC^{m,\a}$ and $\cC_{(\g,H)}^{m,\a}$, where the latter 
is the subset of $\cC^{m,\a}$ consisting of pairs $(g,u)$ where $g$ induces Bartnik boundary data $(\g, H)$.

\begin{proposition} \label{submers}
The scalar curvature map 
$$s: \cS^{m,\a}(M) \to C_{\d+2}^{m-2,\a}(M), \ \ (g, u) \mapsto s(g)$$
is a smooth submersion at any $(g, u) \in \cS^{m,\a}(M)$, i.e.~the linearization $Ds_{g}$ is surjective and its kernel splits. 
The same statement holds for the restricted map 
\be \label{ssub}
s: \cS_{(\g,H)}^{m,\a}(M) \to C_{\d+2}^{m-2,\a}(M), \ \ (g, u) \mapsto s(g).
\ee

\noindent Consequently, the spaces $\cC^{m,\a}$ and $\cC_{(\g,H)}^{m,\a}$ (the latter if nonempty) are smooth Banach manifolds, (closed submanifolds 
of $\cS^{m,\a}(M)$). 
\end{proposition}
  
   A similar result was proved by Bartnik \cite[Theorem 3.7]{Ba4} for complete AF manifolds in a Hilbert space setting for the 
general constraint equations. The proof below is conceptually related. On the one hand, it is simpler than Bartnik's since one only 
has to take account of the scalar constraint ($s = 0$); on the other hand it is more difficult, due to the presence of 
boundary conditions. 

\smallskip 
\noindent {\bf Proof:}  We prove the second statement (i.e., for $\cS_{(\g,H)}^{m,\a}(M)$), which implies the first (for $\cS^{m,\a}(M)$). 
Note also that both statements are independent of $u$, so we may assume $u>0$.

  The proof proceeds (of course) by the implicit function theorem in Banach spaces. To apply this, one needs to prove that the 
linearization $s'=Ds_g$ in \eqref{ssub}, i.e. 
\be \label{cr}
s': T\cS_{(\g,H)}^{m,\a}(M) \to TC_{\d+2}^{m-2,\a}(M), \ \ (h, u') \to s'(h),
\ee
is surjective and the kernel $\Ker s'$ splits as a subspace of $T\cS_{(\g,H)}^{m,\a}(M)$ at $(g,u)$. The usual proofs of these  
properties for compact manifolds or complete AF manifolds, based either on conformal deformations, or on the structure of the 
formally elliptic operator $s'\circ (s')^{*}$, will not work in this setting due to the presence of the boundary conditions. 
We begin with the proof of surjectivity. 
   
   Consider the $4$-manifold $\cM$ and static metrics $g_{\cM}$ on $\cM$, as in \eqref{N}. The Ricci curvature of the 
metric $g_{\cM}$, determined by $(g, u)$, is given by 
$$\Ric_{g_{\cM}} = \Ric_{g} - u^{-1}D^{2}u + u^{-1}\D u \omega \otimes \omega,$$
where $\omega = u dt$ is a $g_{\cM}$-unit covector in the ``vertical'' direction and the first two terms are in the ``horizontal" direction (tangent to $M$).
Passing to the associated Einstein tensor $\Ric_{g_{\cM}} - \frac{s_{g_{\cM}}}{2}g_{\cM}$ 
of $g_{\cM}$ gives 
\be \label{E}
E_{g_{\cM}} = E_{g} - u^{-1}(D^{2}u - \D u \,g) + \frac{s}{2} \omega \otimes \omega,
\ee
where $E_g$ is the Einstein tensor of $g$. 
Consider now the divergence-gauged Einstein operator at a background metric $\w g_{\cM} \in \cS^{m,\a}(M)$:  
$$\Phi_{\w g_{\cM}}: \cS_{\d}^{m,\a}(M) \to \cS_{\d+2}^{m-2,\a}(M),$$
\be \label{Phi}
\Phi_{\w g_{\cM}}(g_{\cM}) = E_{g_{\cM}} + 2\d_{g_{\cM}}^{*}\d_{\w g_{\cM}}(g_{\cM}),
\ee
where $\d_{\w g_{\cM}}$ is the divergence operator with respect to $\w g_{\cM}$ and $\d_{g_{\cM}}^{*}$ is the formal 
$L^2$-adjoint of $\d_{g_{\cM}}$. (Here, we view $ \cS^{m-2,\a}_{\d+2}(M)$ as the Banach space of 
pairs $(\tau,f)$ 
corresponding to $\tau + f \omega \otimes \omega$, 
where $\tau$ is a symmetric 2-tensor on $M$ and $f$ is a function on $M$ with $\tau_{ij}$ and $f$ in 
$C^{m-2,\a}_{\d+2}(M)$). Let $L = D_{\w g_{\cM}}\Phi$ be linearization of $\Phi := \Phi_{\w g_{\cM}}$ at $g_{\cM} = \w g_{\cM}$ and 
let $\hat h = (h, u')$ denote a variation of $g_{\cM}$ with $h$ the variation of $g$ and $u'$ the variation of $u$. We then set 
$$\w L(h, u') = L(h, u') - \D_{g}\nu' \omega \otimes \omega,$$
where $\nu' = (\log u)' = \frac{u'}{u}$. Thus the vertical component of $\w L$ is given by $\frac{1}{2}s'(h) - \D_{g}\nu' + 
\d_{g_{\cM}}^{*}\d_{g_{\cM}} h(V,V)$, where $V = u^{-1}\partial_{t}$, (compare with the remark following \eqref{cL}). We claim 
that $\w L$ is a self-adjoint elliptic operator with respect to the boundary conditions 
\be \label{bc}
(\d_{g_{\cM}}\hat h, h^{T}, H'_{h} + N(\nu')) = (0,0,0) \ \ {\rm at} \ \ \dm. 
\ee
(We note that $H'_{h} + N(\nu')$ is the variation of the mean curvature of $\partial \cM$ in the direction $\hat h$). 
Namely, it is proved in \cite[Lemma 3.2]{AK} that the operator $L$ is self-adjoint and elliptic with respect to the closely 
related boundary conditions 
\be \label{bc2}
(\d_{g_{\cM}}\hat h, [\hat h^{T}]_{0}, H'_{h} + N(\nu')) = (0,0,0) \ \ {\rm at} \ \ \dm,
\ee
where $[\hat h^{T}]_{0}$ denotes the trace-free part of $\hat h^{T}$ (where $\hat h^{T}$ is the restriction of $\hat h$ to $\partial \cM$, and the trace is taken with respect to the induced metric on $\partial \cM$). Using exactly the same methods as in 
\cite[Lemma 3.2]{AK} together with \cite[Proposition 3.7]{AK}, one easily sees that adding the term $-\D_{g}\nu'\omega \otimes \omega$ to 
$L$ has the effect of changing the boundary conditions \eqref{bc2} to \eqref{bc} while preserving ellipticity and self-adjointness. 

  Let $T_{0} \subset T(\cS_{\d}^{m,\a}(M))$ be the subspace of tangent vectors (infinitesimal variations) satisfying the boundary 
conditions \eqref{bc}. By elliptic theory, the operator $\w L_{0} := \w L|_{T_{0}}: T_{0} \to T(\cS_{\d+2}^{m-2,\a}(M))$ is 
Fredholm, so has closed range, and 
\be \label{sa}
({\rm Im} \, \w L_0)^{\perp} = {\rm Ker} \, \w L_0. 
\ee
Let $\pi_{2}$ be the projection onto the $2^{\rm nd}$ (i.e.~vertical) factor of $T(\cS_{\d+2}^{m-2,\a}(M))$, as in \eqref{E}. 
A simple calculation (cf.~\cite[(2.16)]{AK} for instance) shows that $2\d_{g_{\cM}}^{*}\d_{g_{\cM}}h(V, V) = 
\d h(\nabla \log u)$. Hence, 
\be \label{vert}
(\pi_{2}\circ \w L_{0})(h, u') = {\tfrac{1}{2}}s'(h) - \D \nu' + {\tfrac{1}{2}}\d h(\nabla \log u).
\ee
As noted in the beginning of the proof, the statement \eqref{ssub} does not depend on the choice of $u$. Hence, we may set 
$u = 1$, so that $(\pi_{2}\circ\w L_{0})(h, u') = \frac{1}{2}s'(h) - \D u'$. 

  Now the surjectivity of $s'$ on $T\cS_{(\g,H)}^{m,\a}(M)$ is equivalent to the statement that the restriction of $\w L_{0}$ to 
the horizontal subspace $(h, 0)$ (where $u' = 0$) of $T_{0}$ is surjective onto the vertical space $(0,f)$. By the self-adjoint 
property \eqref{sa}, this is equivalent to showing that if 
\be \label{zphi0}
\int_{M}\f s'(h)=0,
\ee
for some $\f \in C^{m,\a}_{\d}(M)$ and for all $(h,0) \in T_{0}$, then $\f = 0$. To prove this, it follows from \eqref{zphi0} that for all such $h$,
\be \label{zphi}
0 = \int_M \f s'(h) = \int_M \<S^{*}\f, h\>,
\ee
where, as in \eqref{S*}, $S^{*}$ is the formal $L^{2}$-adjoint of the linearization $s'$. (The condition $(h,0) \in T_{0}$ 
implies the boundary terms in \eqref{zphi} vanish at $\dm$ and at infinity). In particular $(g,\f)$ solves 
$$S^* \f = 0$$
on $M$ so that $(g, \f)$ is a static vacuum solution, (cf.~the discussion following \eqref{stat1}). 

  Now the potential $\f$ decays to zero at infinity at a rate at least $r^{-\d}$. A simple asymptotic analysis of the 
static vacuum equation \eqref{S*} then implies $\f = 0$ on $M$. Alternatively, this follows directly from Proposition 2.1 of 
\cite{BS} which shows that such potential functions, if not identically zero, are necessarily asymptotic to a non-zero constant or an 
affine function. This completes the proof that $s'$ is surjective. 

   Finally we prove that the kernel $\Ker s'$ splits. Consider again the horizontal and vertical decomposition of the target space 
$\cS_{\d+2}^{m-2,\a}(M)$ as in \eqref{E}. Define $S_{1} = (\w L_0)^{-1}(*, 0)$ and $S_{2} = (\w L_0)^{-1}(0,*)$. Clearly 
$S_{1}$ and $S_{2}$ are closed subspaces of the domain $T_{0}$ of $\w L_{0}$. It is easy to see that $S_1 + S_2$ is also a closed subspace, 
of finite codimension (the latter since the range of $\w L_0$ has finite codimension). Thus $S_1+S_2$ admits a closed complement, 
$S_3$:
$$T_{0} = (S_{1} + S_{2})\oplus S_3.$$
Now, the intersection $S_{1} \cap S_{2}$ equals $\Ker \w L_0$, which is finite dimensional since $\w L_0$ is Fredholm. Hence 
$\Ker \w L_0 \subset S_{2}$ has a closed complement $S_{2}'$ in $S_{2}$. This gives a direct sum 
decomposition 
\be \label{split}
T_{0} = S_{1} \oplus S_{2}' \oplus S_3.
\ee
As in \eqref{vert}, $S_{1} = \{\hat h=(h,u'): \frac{1}{2}s'(h) - \D u' = 0\}$ when $u = 1$. 
Thus, $\Ker s' \cap T_{0}$ is the horizontal subspace of $S_{1}$, which splits by its vertical complement in $S_{1}$. This proves that 
$\Ker s' \cap T_{0}$ splits in $T_{0}$. To complete the proof, let $\chi^{m-1,\a}(M)$ be the space of $C^{m-1,\a}$ vector fields on $M$ 
and let $\chi_{\dm}^{m-1,\a}(M)$ be the restriction of such vector fields to $\dm$. Choose a smooth extension operator $E: \chi_{\dm}^{m-1,\a}(M) 
\to \chi^{m-1,\a}(M)$ and a $C^{m+1,\a}$ smooth metric $\w g$ on $M$ near $(M, g)$. Given $h$, consider the vector field 
$Y_{0}$ dual to $\d h|_{\dm}$ along $\dm$ and let $Y = E(Y_{0})$. One may then uniquely solve the equation $\d_{g}\d_{\w g}^{*}X = Y$ with 
$X = 0$ on $\dm$ and the decomposition $h = (h - \d_{\w g}^{*}X) + \d_{\w g}^{*}X$ gives a 
splitting of $T_{0}$ in the full tangent space $T(\cS_{(\g,H)}^{m,\a}(M))$. This proves that $\Ker s'$ splits.  

 The implicit function theorem (regular value theorem) for Banach manifolds then implies the remaining part of the 
Proposition. 

{\endproof}

\begin{remark}
{\rm  Proposition \ref{submers} is closely related to the issue of linearization stability of solutions of the vacuum Einstein 
equations on $\cM$ and to the work of Fischer, Marsden, and Moncrief \cite{FMM}. While Proposition \ref{submers} is known 
to be false for compact manifolds $M$ (i.e., there exist compact $(M,g)$ for which $Ds_g$ is not surjective), it is known to be true in the $C^{\infty}$ setting for complete AF manifolds (in both 
cases without boundary).

}
\end{remark}

  Proposition \ref{submers} has the following useful corollary. Consider the map to Bartnik boundary data: 
\be \label{pib}
\pi_{B}: \cC^{m,\a}(M) \to Met^{m,\a}(S^{2})\times C^{m-1,\a}(S^{2}), \ \ \pi_{B}(g,u) = (\g, H),
\ee
where $Met^{m,\a}(S^{2})$ is the space of $C^{m,\a}$ Riemannian metrics on $S^2$ with the $C^{m,\a}$ topology. 
Proposition \ref{submers} shows that $\cC^{m,\a}(M)$ is a smooth Banach manifold; clearly $\pi_{B}$ is a smooth 
map of Banach manifolds. 

\begin{corollary} \label{sub}
The map $\pi_{B}$ in \eqref{pib} is a submersion, i.e.~$D\pi_{B}$ is surjective, with splitting kernel. In particular $\pi_{B}$ 
is an open map. 
\end{corollary}

\noindent {\bf Proof:} Proposition \ref{submers} implies that for any $(\g, H) \in \Image \pi_{B}$, the map $s$ is a submersion on 
$\cS_{(\g, H)}^{m,\a}(M)$. Thus, for any given $g$ (or, more precisely, $(g, u)$)  in  $\pi_{B}^{-1}(\g, H)$ and for any $f \in C_{\d+2}^{m-2,\a}(M)$ 
there exists an $h$ (or $(h, u')$) in $T_{(g,u)}\cS^{m,\a}(M)$ satisfying $(h^{T}, H'_{h}) = (0,0)$ at $\dm$ and such that 
$s'(h) = f$. Now given an arbitrary boundary variation $(h^{T}, H'_{h})$, let $h_{e} \in T_{(g,u)}\cS^{m,\a}(M)$ be an extension 
of $(h^{T}, H'_{h})$ to a variation of $g$ on $M$ of compact support. Then $s'(h_{e}) = \f_{e}$, for some $\f_{e} \in C_{\d+2}^{m-2,\a}(M)$. 
Let $h_{0}$ be a solution of $s'(h_{0}) = \f_{e}$ with zero boundary data. Then $h := h_{e} - h_{0}$ satisfies $s'(h) = 0$ 
and $h$ has the given boundary data $(h^{T}, H'_{h})$. This proves that $D\pi_{B}$ is surjective. 

  Further, one has $\Ker D\pi_{B} = T\cC_{(\g,H)}^{m,\a}$ which was proved to split in Proposition \ref{submers}. 
  
{\endproof}

  More generally, for a given  function $\s \in C_{\d+2}^{m-2,\a}(M)$, let 
\be \label{pib1}
\cC_{\s}^{m,\a}(M) = \{g \in Met_{\d}^{m,\a}(M): s_{g} = \s\}.
\ee
One has a corresponding map $\pi_{B}$ as in \eqref{pib} and the same proof as above shows that $\pi_{B}$ remains a submersion, 
for any $\s$. 

\medskip 

  Although the mass $m_{ADM}$ has no critical points on $\cS_{(\g,H)}^{m,\a}(M)$, it may have critical points on distinguished 
subsets of $\cS_{(\g,H)}^{m,\a}(M)$ (constrained critical points). In view of Lemma \ref{lemma_grad_m}, we focus in particular on the submanifold 
$\cC_{(\g, H)}^{m,\a}$. Clearly 
$$m_{ADM}: \cC_{(\g,H)}^{m,\a} \to \bR$$
is a smooth functional. 

\begin{proposition}
\label{prop_critical}
Critical points of the ADM mass $m_{ADM}$ on $\cC_{(\g,H)}^{m,\a}$ are exactly metrics $g$ that admit an AF function $u$ such that 
$(g,u)$ is an AF solution of the static vacuum Einstein equations on $M$ with given boundary data $(\g, H)$. 
\end{proposition}
Note that we are not yet claiming that $u>0$; this is addressed in Theorem \ref{thm_BC} for Bartnik mass minimizers.

\noindent {\bf Proof:} On $\cC_{(\g,H)}^{m,\a}$ (if nonempty), one has by \eqref{cL} 
\be \label{Hem}
\cH = -16\pi m_{ADM}: \cC_{(\g,H)}^{m,\a} \to \bR.
\ee
The critical points of $m_{ADM}$ on the constraint space $\cC_{(\g,H)}^{m,\a}$ are thus exactly the same as 
critical points of $\cH$ on $\cC_{(\g,H)}^{m,\a}$, and in the following we work with $\cH$. Note that the potential 
function $u$ is irrelevant at this point (since $m_{ADM}$ is independent of $u$); thus in the following and for the 
moment, we make a fixed (but arbitrary) choice of $u = u_{0}>0$, with $u_{0} - 1 \in C_{\d}^{m,\a}(M)$. 

  Let then $(g, u_{0})$ be a critical point of the constrained variational problem, i.e.
$$d\cH_{(g,u_{0})}(h, u') = 0,$$
for all $(h,u') \in T_{(g,u_{0})}\cC_{(\g,H)}^{m,\a}$, i.e. $s'(h) = 0$ and $(h^{T}, H'_{h}) = (0,0)$ at $\dm$. 
A standard Lagrange multiplier theorem, discussed explicitly in a related context in \cite[Theorem 6.3]{Ba4}, shows that there is 
a distribution (bounded linear functional) $\l$ on $C_{\d+2}^{m-2,\a}(M)$, (the Lagrange multiplier), such that for all, i.e.~unconstrained, variations 
$(h, u') \in T_{(g,u_{0})} \cS^{m,\a}_{(\g,H)}$, 
$$d\cH_{(g, u_{0})}(h, u') = \l(s'(h)).$$
Since $s_{g} = 0$, one has by Proposition \ref{l2grad},
$$d\cH_{(g,u_{0})}(h,u') = \int_M \langle S^*u_0, h \rangle = \int_M u_0 s'(h),$$ 
for all variations $(h,u') \in T\cS_{(\g,H)}^{m,\a}$ of compact support in $\interior(M)$. Combining these statements gives 
\be\label{lambda_u_0}
\l(s'(h))-\int_M u_0 s'(h) =0,
\ee
for all such variations $h$. We claim that the distribution 
$$T(w):=\int_Mu_0 w - \lambda(w),$$
acting on functions $w \in C^{m-2,\a}(M)$ with compact support in $\interior(M)$ is a smooth solution of the static vacuum equations with respect to $g$ on $M$, $C^{m,\a}$ up to $\dm$. 
To see this,  consider variations of the form $h = fg$, where $f$ is a smooth function on $M$. Since $s'(fg) = -2\Delta f -s f=-2\Delta f$, 
we see from \eqref{lambda_u_0} that $T$ is a weak (i.e.~distribution) solution of $-2\Delta u = 0$ in $\interior(M)$. By elliptic regularity (i.e.~the well-known 
Weyl Lemma and Schauder estimates), $T$ is represented as
a $C^{m,\alpha}$ function $u$, i.e.,
$$T(w) = \int_M u w$$
for all such $w$, where $\Delta u = 0$. Thus for all $ (h, u') \in T_{(g,u_{0})} \cS^{m,\a}_{(\g,H)}$ 
of compact support in $\interior(M)$, 
$$0 = \int_M u s'(h) = \int_M \<S^{*}u, h\>,$$
where, as in \eqref{S*}, $S^{*}$ is the formal $L^{2}$-adjoint of the linearization $s'$. 
In particular $(g, u)$ solves $S^* u = 0$ in $\interior(M)$ so that $(g, u)$ is a smooth static vacuum solution, (cf.~the discussion 
following \eqref{stat1}). Since $\Ric_{g} \in C^{m-2,\a}(M)$, integration of the static equations $u\Ric = D^{2}u$ shows that 
$u$ is $C^{m,\a}$ up to $\dm$.

  We complete the proof by arguing $u - 1 \in C_{\d}^{m,\a}(M)$. Note
$$\lambda(w) = \int_M (u_0 -u)w,$$
for $w$ compactly supported in $\interior(M)$, where we recall $u_0 \to 1$ at infinity. From \cite[Proposition 2.1]{BS}, the static vacuum potential $u$ must approach a nonzero constant at infinity, be asymptotic to an affine function at infinity, or else be identically zero. We will reach a contradiction (to $\lambda$ being a bounded linear functional) in every case except $u \to 1$ at infinity. This together with $\Delta u= 0$ on $M$ implies $u -1 \in C_{\d}^{m,\a}(M)$, and the proof will be complete.

Let $\chi_{i} : M \to [0,1]$ be a sequence of smooth, nonnegative radially symmetric 
cut-off functions of compact support on $M$, with $\chi_{i}(r) = 1$ for $r \in [0, R_{i}]$, $|d^{k}\chi_{i}(r)| \leq c_{k}/r^{k}$ for 
$r \in [R_{i}, 2R_{i}]$ and $\chi_{i}(r) = 0$ for $r \geq 2R_{i}$ with $R_{i} \to \infty$. 

First, suppose that either $u \equiv 0$, or $u$ approaches a constant other than 1 at infinity. Then $u_0 - u$ approaches a nonzero constant at infinity. Let $w \in C_{\d+2}^{m-2,\a}(M)$ equal $r^{-2-\delta}$ outside of a compact set and vanish near $\partial M$.  
Let $w_i = \chi_i w$, a uniformly bounded sequence in $C_{\d+2}^{m-2,\a}(M)$ each with compact support in $\interior(M)$. Since $\l$ is a bounded linear functional, the sequence $\l(w_i)$ is uniformly bounded as well. However, a direct computation shows that $|\l(w_i)|$ diverges to infinity, a contradiction.

Last, consider the case in which $u$ (and hence $u_0-u$) limits to an affine function at infinity, i.e. 
$$(u_0-u)(\vec x) = \vec a \cdot \vec x + O(|\vec x|^{1-\delta}),$$ 
where $\vec a \neq \vec 0$ is constant (cf.~\cite{BS}). Without loss of generality, assume $\vec a = (0,0,1)$, i.e. $(u_0-u)(x,y,z) = z + O(|\vec x|^{1-\delta})$.
A similar argument as that given above applies here, with a different sequence of test functions. Let
$$w_i(\vec x) = |\vec x|^{-2-\d} \chi_i (\vec x - (0,0,3R_i)),$$
a sequence of smooth functions of compact support in $\interior(M)$ vanishing near $\partial M$ (for $i$ large enough). It is  straightforward to check that $\{w_i\}$ is uniformly bounded in $C^{m-2,\a}_{\d+2}(M)$, but that $\lim_{i \to \infty} \l(w_i) = +\infty$. Again, this is a contradiction, since $\l$ is a bounded linear functional. The proof is complete.

{\endproof}

   We note that the potential $u$ of an AF static vacuum metric $(M, g)$ is uniquely determined by $g$ (up to multiplication by a 
scalar) if $g$ is not flat, cf.~\cite[Proposition 10]{T}. On the other hand, any affine function $u$ is the potential of a flat 
static exterior solution $(M, g_{Eucl})$. As discussed in Remark \ref{corv} below, it is not fully known in general whether, if $(M, g)$ is an 
AF static metric, then the potential $u$ must also be AF when a (non-mimimal) boundary is present.

  We now state and prove one of the main results of this section:
\begin{theorem}
\label{thm_mass_min_static}
An AF metric $g \in \accentset{\circ} \cP^{m,\a}(M)$ on $M$ realizing the Bartnik mass \eqref{bm21} of the boundary data $(\g, H)$ is an AF static vacuum 
solution $(g, u)$ (i.e., with $u \to 1$ at infinity), satisfying the boundary conditions \eqref{bcont_weak}. 
\end{theorem}

\noindent {\bf Proof:} The Bartnik mass $m_{B}$ is obtained by minimizing $m_{ADM}$ subject to the no-horizon condition and 
constraints that $s \geq 0$, the boundary metric $\g$ is fixed, and the mean curvature of $\dm$ is at most $H = H_{\partial \O}$ pointwise. 
If $g \in \accentset{\circ} \cP^{m,\a}(M)$ realizes $m_{B}$ (say with mean curvature $H_{\dm}$), then a neighborhood of $g$ in $\cP^{m,\a}(M)$ is in $\accentset{\circ} \cP^{m,\a}(M)$ 
by Lemma \ref{lemma_horizons}. Lemma \ref{lemma_grad_m} then implies that $g$ must be scalar-flat, $s_{g} = 0$. The result 
then follows from Proposition \ref{prop_critical}, since  $g$ is a critical point of  $m_{ADM}$ on $\cC_{(\g,H_{\dm})}^{m,\a}$.

{\endproof}

\begin{remark}
\label{remark_BCs}
{\rm We recall briefly the reasoning that leads to the Bartnik boundary conditions \eqref{bcont} and \eqref{bcont_weak}. 
By combining the Gauss and Ricatti equations on $M$ at $\dm = \partial \O$ one finds 
\be \label{nh}
N(H) = {\tfrac{1}{2}}(s_{\g} - s_{g} - |A|^{2} - H^{2}).
\ee
Since the metric $\g$ is fixed on $\partial \Omega$, the scalar curvature $s_{\g}$ is fixed, while the last three terms in 
\eqref{nh} are non-positive, since $s_{g} \geq 0$. It follows that $N(H)$ is uniformly bounded above, but may a priori become arbitrarily 
negative in (weak) limits. Thus in passing to a limit of a mass-minimizing sequence of extensions, one expects 
\be \label{Hless}
H_{\partial M} \leq H_{\partial \O},
\ee
so that the exterior mean curvature may drop from that given by the region $\O$, as in \eqref{bcont_weak}. Note that the 
positive mass theorem still holds on such manifolds with corners, cf.~\cite{Mi1}, \cite{ST}. On the other hand, if $A$ 
and $s_{g}$ remain bounded on a minimizing sequence, then one has 
\be \label{equal}
H_{\partial M} = H_{\partial \O}.
\ee
Unfortunately, it is not clear how to give a topology on $\cP^{m,\a}(M)$ to effectively 
implement such a structure in limits. 

  Conversely, given (say) smooth boundary data $(\g, H)$ on $S^{2}$ that arise as boundary data for a smooth 
metric $g_{\Omega}$ of non-negative scalar curvature on $\Omega \simeq \bar B$, one expects that there is a smooth 
AF metric on $\bR^{3}$ of non-negative scalar curvature in which $(\Omega, g_{\Omega})$ isometrically embeds, 
(corresponding to the original definition \eqref{bm}). This remains to be fully proved however. 

} 
\end{remark}

   Next we prove that boundary conditions \eqref{equal} are actually realized for a minimizer of the ADM mass, as 
defined in \eqref{bcont_weak}, provided one uses the definition \eqref{P0} for $\accentset{\circ} \cP^{m,\a}(M)$. This result was obtained 
by Miao for the case in which $(\gamma, H)$ has strictly positive Gauss curvature (and $H>0$) using a different technique; 
see \cite[Proposition 3.4]{M1}. In addition, we prove that $u > 0$ everywhere on $M$. 
   
\begin{theorem}
\label{thm_BC}
Suppose that an AF metric $g$ on $M$ realizes the Bartnik mass of the boundary data $(\g, H)$  in the sense of \eqref{bm21} (i.e.~with 
boundary conditions \eqref{Hless}), where $H>0$. Then \eqref{equal} holds, and the AF static vacuum potential $u$ is strictly positive on $M$.
\end{theorem}

\noindent {\bf Proof:}  Suppose $(M, g)$ realizes the Bartnik mass \eqref{bm21} of $(\O, g_{\O})$, so that as in \eqref{Hless}, 
$H_{\partial M} \leq H_{\partial \O}$. We first show that \eqref{equal} holds. If it fails let $U \subset \dm$ be the nonempty open set on which strict inequality holds:
$$H_{\partial M} < H_{\partial \O} \ \ \text{ on } U.$$
By Theorem \ref{thm_mass_min_static}, the metric $g$ is static vacuum, and so in particular scalar-flat. Let $u$ be the corresponding AF
static vacuum potential.

   Consider the map $\pi_{B}$ in \eqref{pib}. By Corollary \ref{sub}, $\pi_{B}$ is a submersion at $(g, 1)$ and so for any 
variation $(h^{T}, H'_{h})$ of the boundary data $(\g, H_{\partial M})$, there is a variation $h$ of $g$ such that $s'(h) = 0$. 
Choose 
\be \label{q}
(h^{T}, H'_{h}) = (0, q),
\ee
where $q$ is a smooth function on $\dm$ supported in $U$ and such that 
\be \label{q2}
\int_{\dm}u q  > 0.
\ee
Clearly, there are many such choices of $q$, unless $u \equiv 0$ on $U$. However, if $u \equiv 0$ on $U$, then $H_{\partial M} \equiv 0$ on $U$, 
since the zero set of a static vacuum potential is totally geodesic. Consider two cases. First, if $U$ is a proper subset of $\partial M$, then 
$H_{\partial M} \equiv 0$ on $U$ contradicts $H_{\partial \Omega} > 0$, since $H_{\partial M} = H_{\partial \O}$ outside $U$. Second, if $U= \dm$, then 
$(M, g, u)$ is a Schwarzschild metric, by the black hole uniqueness theorem \cite[Theorem 2]{BM}. In particular, $\g = \g_{2m}$ is a round metric. 
Since $H_{\partial \Omega}>0$ and $H_{\partial M}=0$, it is easy to see that $(M, g)$ cannot be a minimal mass extension. 
(For example, one may take an equidistant round sphere $r > 2m$ close to the horizon $r = 2m$ of the Schwarzschild metric 
and rescale, decreasing the mass). 

Now, let $h$ be the corresponding variation of $g$ with $s'(h)=0$, satisfying \eqref{q}. Since $s'(h) = 0$ and $(g,u)$ is static vacuum,
one has from \eqref{Hem}, \eqref{var}, and \eqref{q} that 
\be \label{ml00}
-16\pi m_{ADM}'(h) = d\cH_{(g,u)}(h,0) = 2\int_{\dm}uq > 0,
\ee
so that $m'_{ADM}(h) < 0$. This gives, at the infinitesimal level, a mass-decreasing variation of $(g, 1)$ in 
$T_{(g, 1)}\cC_{(\g, H_{\partial \O}^{\leq})}^{m,\a}(M)$. Now consider the curve of boundary data $(\g, H_{\partial M} + tq)$ 
(for instance) with $t$ small. Again by Corollary \ref{sub}, this curve lifts (via $\pi_{B}^{-1}$) to a curve $g_{t}$ 
in the slice or closed complement to $\Ker s' = T_{(g,1)}\cC^{m,\a}$ in $\cC^{m,\a}$. It follows that 
$m_{ADM}(g_{t}) < m_{ADM}(g)$ for $t > 0$ small and since $\pi_{B}(g_{t}) = (\g, H+tq)$, $g_{t} \in \cC_{(\g, H_{\partial \O}^{\leq})}^{m,\a}(M)$, 
again for $t$ small. Since by Lemma \ref{lemma_horizons}, $g_{t}$ has no horizons for $t$ small, this contradicts the definition of Bartnik mass. Thus \eqref{equal} holds.

Next, suppose $u < 0$ somewhere on $M$. By the maximum principle, $u<0$ at some point on $\dm$. Let $q$ be a smooth, non-positive function on $\dm$, supported in the set where $u<0$, satisfying \eqref{q2}. A similar argument to that given above produces a mass-decreasing path of metrics in $\cC_{(\g, H_{\partial \O}^{\leq})}^{m,\a}(M),$ again contradicting the definition of the Bartnik mass. Thus $u \geq 0$ on $M$.

Finally, if $u(p)=0$ for some $p \in M$, then $p \in \dm$ by the maximum principle. At $p$, by the static vacuum equations, $0 = u\Ric = D^{2}u$ and hence $D^{2}u = 0$ at $p$. The restriction of $D^{2}u$ to $\dm$ gives $(D^{2})^{T}u + N(u)A = 0$ and taking then the trace over $\dm$ gives $\D_{\dm}u + N(u)H = 0$. Since $p$ is a minimum of 
$u$, $\D_{\dm}u \geq 0$, while by the Hopf boundary point maximum principle, $N(u) > 0$. It follows then that 
$H(p) \leq 0$, a contradiction. This proves $u>0$.
 
{\endproof} 

Theorems \ref{thm_mass_min_static} and \ref{thm_BC} together imply Theorem \ref{thm1.1} from the introduction.

\begin{remark} \label{corv}
{\rm  Theorem \ref{thm1.1} and its proof implement the original heuristic program suggested by Bartnik \cite{Ba2}, \cite{Ba4}, 
based on the perspective initiated by Brill--Deser--Fadeev \cite{BDF}. It is also possible to prove part of Theorem \ref{thm1.1} using different methods. We briefly summarize this alternate approach and some related results below.

To begin, Corvino \cite{C1} showed that metrics minimizing the Bartnik mass of a domain $\O$ are static vacuum outside 
$\bar \O$ by constructing suitable localized scalar curvature deformations. However, this result did not address the issues 
of the horizon conditions, nor the global behavior of the potential function $u$, and did not fully address the boundary conditions.

  Using this method, an elementary argument in \cite{C1} shows that the boundary condition \eqref{bcont_weak} is preserved; however, 
it is not clear whether the original (stronger) condition \eqref{bcont} is preserved. This issue has very recently been resolved in \cite{C2}. 
Note that the proof that a minimizer is static vacuum does require some stability condition, such as that in Lemma \ref{lemma_horizons}. 

  However, in such an approach it is not clear if the potential function $u$ satisfies $u > 0$ or even $u \to 1$ at infinity $M$, i.e.~
$(g, u)$ may not be globally static vacuum or asymptotically flat (AF) in the usual sense. Miao--Tam \cite[Theorem 1.1]{Mi3} prove that a static potential is bounded, with a sign, outside a compact set, provided the metric is asymptotically Schwarzschild of nonzero mass. Unfortunately, this does not apply to the present setting, since it is not clear that a static vacuum metric is asymptotically Schwarzschild without knowing in advance its potential has a sign at infinity. Huang--Martin--Miao \cite[Theorem 10]{HMM} show that the static potential of a Bartnik mass minimizer approaches 1 at infinity (though this argument does not show $u>0$ globally). Note that they use a different version of the Bartnik mass. We also refer the reader to some related results on static potentials by Galloway--Miao \cite{GM} (cf. Carlotto--Chodosh--Eichmair \cite[Corollary 1.9]{CCE}), for example, in the boundaryless or minimal boundary case.}
\end{remark}

  An immediate corollary of Theorem \ref{thm_BC} is the strict monotonicity of the Bartnik mass, improving \eqref{mono} when 
$m_{B}$ is realized. 
  
\begin{corollary} \label{mincor}
Suppose $0 < H \leq H'$ and $H < H'$ on some nonempty open set $U \subset \dm$. If the data $(\g, H)$ is realized by a 
mass-minimizing extension in $\accentset{\circ} \cP^{m,\a}(M)$, then  
\be \label{mono1}
m_{B}(\g, H') < m_{B}(\g, H).
\ee

\end{corollary}

\section{Remarks on Conjecture I}
\setcounter{equation}{0}

  In this section, we present several results that provide further positive evidence for the validity of Conjecture I, regarding the existence 
of AF extensions of non-negative scalar curvature.

  The main extension results to date are based on the quasi-spherical method introduced by Bartnik \cite{Ba5}. For example, using this 
method, it can be established that any boundary data $(\g, H)$ with $\g$ of positive Gauss curvature $K_{\g} > 0$ and $H > 0$ admits 
an extension in $\cP(M)$ (see \cite{Ba5}, \cite{ST}, \cite{SW}). More recently, extension results have also been obtained by 
Lin \cite{Lin} and Lin--Sormani \cite{LS} using a modified Ricci flow. 

  We write $(\g, H) \in \cP^{m,\a}(M)$ if the boundary data $(\g, H)$ on $S^{2}$ admit an admissible extension $g \in 
\cP^{m,\a}(M)$, and similarly for $\accentset{\circ} \cP^{m,\a}(M)$. 

  We first note the following general result. 
\begin{proposition} \label{sopen}
The spaces $\cP^{m,\a}(M)$ and $\accentset{\circ} \cP^{m,\a}(M)$ are open in $Met^{m,\a}(S^{2}) \times C^{m-1,\a}(S^{2})$. 
\end{proposition}

\noindent {\bf Proof:} This is an immediate consequence of Corollary \ref{sub}, in the scalar-flat case. The general case follows as in 
\eqref{pib1}. Lemma \ref{lemma_horizons} then implies the statement for $\accentset{\circ} \cP^{m,\a}(M)$. 

{\endproof}

   For the discussion to follow, we will not address the horizon issue, which is more difficult to understand when dealing with 
more global problems.  

\medskip 

  We first prove a general result that the space $\cP^{m,\a}(M)$ is invariant under pointwise increase of the mean curvature $H$, keeping 
the boundary metric $\g$ fixed, (compare with the proof of Theorem \ref{thm_BC}). In fact, even a small decrease on $H$ is allowed. 
The method of proof will also be used in the proof of Theorem \ref{thm_conj_II} given in Section 4.

\begin{proposition}
\label{prop_H}
Suppose $(\g, H) \in \cP^{m,\a}(M)$. There exists a $C^{m-1,\a}$ function $\mu>0$ on $S^2$, (depending on $(\g, H)$), 
such that for any $C^{m-1,\a}$ function $H_0$ on $S^2$ satisfying
\be \label{bigger}
H_0 \geq H-\mu,
\ee
pointwise, one has $(\g, H_0) \in \cP^{m,\a}(M)$.
\end{proposition}

\noindent {\bf Proof:}  Suppose $(\g, H) \in \cP^{m,\a}(M)$, so that there is an AF extension $g$ of $(\g, H)$ with scalar curvature 
$s \geq 0$. For simplicity, we will assume $(\g, H)$ and $g$ are smooth. We first consider the case $H_0 \geq H$ and construct 
an AF extension of $(\g, H_0)$ by a conformal deformation of $g$.

  For a conformal metric $\w g = v^{4}g$ with $v>0$, the scalar curvature $s$ of $g$ changes as 
\be \label{confs}
v^{5}\w s = - 8\D v + sv =: f,
\ee
where $\D$ is the Laplacian operator on $(M,g)$. Since $s \geq 0$, $- 8\D + s$ is a positive operator (for Dirichlet boundary data). 
Clearly $\w g \in \cP(M)$ requires $f \geq 0$. 

For $r>0$, let $B(r) = \{x \in M: \dist(x, \dm) \leq r\}$, so that $\partial B(r) = \dm \cup S(r)$ is a regular hypersurface for $r$ large. Given $f\geq 0$ of compact support 
on $M$, let $v_{r}$ be the unique solution to \eqref{confs} on $B(r)$ with Dirichlet boundary data $v_r = 1$ on $\dm \cup S(r)$. 
By the maximum principle, $v_r > 0$ on $B(r)$. It is standard that, letting $r \to \infty$, $v_{r} \to v$ with $v > 0$ on $M$, 
$v = 1$ at $\dm$, and $v \to 1$ at infinity. In particular, $\w g$ is a conformal AF metric on $M$, with induced boundary metric
$$\w \g = \g \ \ {\rm at} \ \  \dm,$$ 
and boundary mean curvature
$$\w H = v^{-2}H + 4v^{-2}N(\log v),$$
where $N$ is the unit normal into $M$ with respect to $g$. Hence  
\be \label{wHH}
\w H = H + 4N(v).
\ee
To obtain $\w H = H_0$, we will demonstrate how to choose $f \geq 0$ appropriately so that the solution $v$ to \eqref{confs} with boundary conditions of 1 
on $\dm$ and at infinity satisfies $N(v) = \frac{1}{4}(H_0-H)$.

Write \eqref{confs} in the form 
\be \label{L}
L(v) := \D v - {\tfrac{1}{8}}sv = -{\tfrac{1}{8}}f.
\ee
On the bounded domain $(B(r), g)$, $L$ has a (negative) Green's function $G$, with $G(x, y) = 0$ for (say) $y \in \dm \cup S(r)$. 
It is standard that we may take $r \to \infty$ to obtain a (negative) Green's function $G$ of $(M,g)$ with $O_2(1/|y|)$ decay at 
infinity for $G(x_0,y)$, for any fixed $x_0$.

The Poisson kernel $P$ of $L$ is given by $P(x, y) = -N_xG(x,y)$, for $x \in \dm$ and is positive for $y$ in the interior of $M$.
Green's formula gives for $x \in M$,
\be \label{rep}
v(x) =  {\tfrac{1}{8}}\int_{M}G(x,y)(-f(y))dy + \int_{\dm}P(y,x)v(y)dy + \lim_{r \to \infty}  \int_{S(r)} N_y G(x,y) v(y)dy.
\ee
Here $dy$ represents the corresponding volume forms on $M$, $\partial M$, and $S(r)$, respectively, with respect to $g$.

Note that $v=1$ solves \eqref{L} uniquely for the choice $f=s$; using \eqref{rep} on the sequence $v_i$ corresponding to $f_i = s \chi_i$, where $0 \leq \chi_i \leq 1$ are compactly supported functions on $M$ that converge pointwise to 1 as $i \to \infty$, we may take a limit (using the decay $O(r^{-2-\delta})$ of $s$) to arrive at
\be \label{one}
1 = {\tfrac{1}{8}}\int_{M}G(x,y)(-s(y))dy + \int_{\dm }P(y,x)1 dy + \lim_{r \to \infty}  \int_{S(r)} N_y G(x,y) 1 dy. 
\ee
Since $v = 1$ on $\dm$ and $v \to 1$ at infinity, it follows that for general $f$ of compact support, the solution $v$ to 
\eqref{L} is given by
\be \label{bulk}
v(x) = 1 + {\tfrac{1}{8}}\int_{M}G(x,y)(s(y)-f(y))dy.
\ee
Thus, for $x \in \dm$, 
\be \label{nv}
N(v)(x) = {\tfrac{1}{8}}\int_{M}P(x,y)(f(y) - s(y))dy.
\ee
It is standard that $v$ has sufficient decay at infinity (e.g., $v(y)-1 = O_m(|y|^{-\delta})$) to assure that $\w g$ is asymptotically flat. 

 Now we claim that given any smooth function $\f \geq 0$ on $\dm$, there is a $C^k$ function $f \geq 0$ (for any $k>0$) with compact 
support on $M$, such that
\be \label{posi}
\f(x) = N(v)(x) = {\tfrac{1}{8}}\int_{M}P(x,y)(f(y) - s(y))dy.
\ee
This will complete the proof (in the case $H_0 \geq H$), by choosing $\f ={\tfrac{1}{4}}( H_0- H) \geq 0$.

To prove the claim, note first that by the basic reproducing property of the Poisson kernel, 
\be \label{repro}
\int_{\dm}P(x,y)\f(y)dy = \f(x).
\ee
Choose a constant $d_{0}>0$, smaller than the distance to the cut-locus of the normal exponential map of $\dm$ 
into $M$, and let $\dm_{r} = \{y \in M: \dist(y, \dm) = r\}$ for $0\leq r \leq d_0$. Define a continuous linear operator 
$A_{r}: C^{k}(\dm_{r}) \to C^{k}(\dm_{r})$ by 
\be \label{surj}
A_{r}(\chi)(x) = \chi_{r}(x) = \int_{\dm_{r}}P_{r}(x,y)\chi(y)dy,
\ee
where $x \in \dm_{r}$ and $P_{r} = P|_{\dm_{r}}$. Using the identification of  $\dm$ with $\dm_{r}$ via the normal exponential map, 
we may regard $A_r$ as a map $C^{k}(\dm) \to C^{k}(\dm)$. It is well-known (and easy to see) that 
$$A_{r} \to Id, \ \ {\rm as} \ \ r \to 0,$$
as bounded linear operators on $C^{k}(\dm)$. Since the space of invertible operators is open, we may shrink $d_0>0$ if necessary 
so that given $\f \in C^{k}(\dm)$ there exists a unique $\f_{r} \in C^{k}(\dm_{r})$, $r \leq d_{0}$, satisfying  
\be \label{f_r}
\int_{\dm_{r}}P_{r}(x,y)\f_{r}(y)dy = \f(x).
\ee
Note that $P_{r}$ is smooth in $r$ for $r > 0$, and hence so is $\f_{r}$. The (higher order) normal derivatives 
$\partial_{r}^{k}P(x,y)$ govern (by convolution) the (higher order) normal derivatives of harmonic functions on 
$M$ at $\dm$; it follows that $\f_{r}$ is also smooth in $r$ at $r = 0$. 

  Now, let $\rho(r)$ be a smooth function of $r \geq 0$ with $\rho(0)=1$, $\rho(r) = 0$ for each $r \geq d_0$, and 
$\int_0^{d_0} \rho(r) = \frac{d_0}{2}$. Integrating over $r$ and using the Gauss Lemma and Fubini theorem 
(or the coarea formula) gives 
\be \label{fubini}
\int_{M}P(x,y)\rho(r(y))\f_{r}(y)dy = \f(x)\int_{0}^{d_{0}}\rho(r)dr = \frac{1}{2}d_{0} \f(x).
\ee
Thus the $C^k$ function given by $f(y) = s(y) + \frac{16}{d_{0}}\rho(r(y))\f_{r}(y)$ satisfies \eqref{posi}.  It is clear 
that $f$ is $C^k$ smooth and extends smoothly by zero to $M$. This proves the claim \eqref{posi}. Note that $f$ is not uniquely 
determined. 
 
 To complete the proof when $\mu > 0$, consider the given extension $(M,g)$ of $(\g, H)$, and take the unique, smooth 
solution to
\[\begin{cases}
\D u = 0 & \text{on } M\\
u =1 & \text{on } \partial M\\
u \to \frac{1}{2} & \text{at } \infty.
\end{cases}\] 
Then the conformal metric $\hat g = u^4 g$ belongs to $\cP^{m,\a}(M)$, induces the metric $\g$ on its boundary, and 
the induced mean curvature on $\partial M$ given by
\[\hat H = H + 4N(u).\]
By the maximum principle, $N(u) < 0$ on $\dm$. Thus,  $(\g, H - \mu) \in \cP^{m,\a}(M)$ for the choice $\mu = -4N(u)>0$. 
The result now follows by applying the argument above to $(\g, H-\mu)$.

{\endproof}

  It is useful to understand how the mass $m_{ADM}$ changes under the deformations in Proposition \ref{prop_H}.  
Thus recall from \eqref{confmass} that if $\w g = v^{4}g$, then 
$$\w m = m - \frac{1}{2\pi}\lim_{r\to \infty}\int_{S(r)}N(v)dA.$$

  In the context of the proof of Proposition \ref{prop_H} above, suppose $f(y) \geq s(y) \geq 0$, so that, 
in particular, $\w H \geq H$. For simplicity, assume $f=s$ outside of a compact set. Then \eqref{bulk} shows $v \geq 1$ and $v \to 1$ at infinity, and 
\eqref{L} shows $v$ is subharmonic outside of a compact set. Then, similar to the proof of the last part of Lemma \ref{lemma_grad_m}, we can enclose any coordinate sphere $S(r)$ with a smooth level set $L$ of $v$, on which the outward normal derivative of $v$ is nonpositive by the maximum principle. Applying the divergence theorem on the region between $S_r$ and $L$, we see that that $\int_{S(r)} N(v) dA \leq 0$ for all $r$ large, and thus
$$\w m \geq m.$$
Thus, roughly speaking, as one increases $H$, the mass $m$ increases under conformal changes, 
when keeping the boundary metric fixed. On other hand, if $0 \leq f(y) \leq s(y)$, so that $\w H \leq H$, 
then \eqref{bulk} shows $v \leq 1$ and $v \to 1$ at infinity, so that 
$$\w m \leq m.$$
In particular, one can decrease the mass conformally if $s \geq 0$ is not identically zero; compare with Lemma \ref{lemma_grad_m}. 

\medskip 

   When combined with existing results, Proposition \ref{prop_H} gives further partial evidence for the validity of Conjecture I. 
  
  Given a metric $\g$ on $S^{2}$, let $\l_{1}(-\D_{\g} + K_{\g}) > 0$ be the lowest eigenvalue of the operator 
$-\D_{\g} + K_{\g}$, where $\D_{\g}$ is the Laplacian with respect to $\g$ and $K_{\g}$ is the Gauss curvature.

\begin{corollary} One has $(\g, H) \in \cP^{m,\a}(M)$ for all $H > 0$ and all $\g$ such that $\l_{1}(-\D + K) > 0$. 
\end{corollary}

\noindent {\bf Proof:} In \cite{MS}, Mantoulidis and Schoen constructed extensions $g \in \cP^{m,\a}(M)$ of $(\g,0)$ for 
$\g$ satisfying $\l_{1}(-\D + K) > 0$. The result then follows from Proposition \ref{prop_H}. 
 
{\endproof}

  This generalizes (with a different proof) previous extension results of Bartnik \cite{Ba5} and Miao \cite{M1}. 

\begin{corollary} 
\label{cor_l_0}
For any $(\g, H)$ with $H > 0$ there is a $\l_{0} > 0$ such that 
$$(\g, \l H) \in \cP^{m,\a}(M),$$
for all $\l \geq \l_{0}$. 
\end{corollary}

\noindent {\bf Proof:}  The proof is based on work in \cite{A1}, \cite{AK} and \cite{A2} on the moduli space $\cE^{m,\a}$ of $C^{m,\a}$ 
AF static vacuum solutions $(g, u)$, $u > 0$, on $M = \bR^{3}\setminus B$. Namely, the map to Bartnik boundary data 
\be \label{PiB10}
\Pi_{B}: \cE^{m, \a} \to Met^{m,\a}(S^{2})\times C^{m-1,\a}(S^{2}),
\ee
$$\Pi_{B}(g,u) = (\g, H),$$
is a smooth Fredholm map, of Fredholm index 0. (This is discussed further in Section 5). Consider the map $\Pi_{B}$ restricted to the 
space $\cE_{+}^{m,\a}$ of static vacuum metrics with $H > 0$ at $\dm$:  
\be \label{Pi+}
\Pi_{B}: \cE_{+}^{m,\a} \to Met^{m,\a}(S^{2})\times C_{+}^{m-1,\a}(M).
\ee
Consider also the action of scalars $\l \in \bR^{+}$ on $C_{+}^{m-1,\a}(S^{2})$ where $(\l, H) \to \l H$. 
Let $\cD_{+}^{m-1,\a}(S^{2})$ be the space of equivalence classes $[H] = [\l H]$. The space 
$\cD_{+}^{m-1,\a}(S^{2})$ is clearly a Banach manifold. 

 It is proved in \cite{A2} that the induced quotient map 
$$\w \Pi_{B}: \cE_{+}^{m, \a} \to Met^{m,\a}(S^{2})\times \cD_{+}^{m-1,\a}(S^{2}),$$
$$\w \Pi_{B}(g) = (\g, [H]),$$ 
is a smooth {\em surjective} Fredholm map of Fredholm index 1. Hence, for any given boundary 
data $(\g, H)$, $H > 0$, there exists $\l_{0} = \l_{0}(\g, H)$ such that $(\g, \l_{0}H)$ are the Bartnik 
boundary data of a complete AF static vacuum solution $(M, g, u)$. Since $g \in \cP^{m,\a}(M)$, 
the result then follows from Proposition \ref{prop_H}. 

{\endproof}

  Corollary \ref{cor_l_0} may be contrasted with the result in \cite{J1} that for given Bartnik boundary data $(\g, H)$, with $H>0$, 
there is a largest value $\l^{0} < \infty$ such that $(\g, \l H)$ has a $s \geq 0$ infilling for $\l < \l_{0}$, and no such infilling, 
for $\l > \l_{0}$. Here, an $s \geq 0$ infilling is a compact Riemannian 3-manifold with boundary inducing Bartnik boundary data 
$(\g, H)$ that has non-negative scalar curvature. That result required $\g$ to have positive Gauss curvature $K_{\g}$, but a recent 
result of Mantoulidis and Miao (\cite[Theorem 1.3]{MM}) implies that $\l^{0} < \infty$ without assuming $K_\gamma >0$.

\section{Proof of Theorem \ref{thm_conj_II}}
\setcounter{equation}{0}

  The main purpose of this section is to prove Theorem \ref{thm_conj_II}. Most of the section will be devoted to proving:

\begin{theorem} 
\label{thm_m0}
Let $F \in \cF$ be as in the statement of Theorem \ref{thm_conj_II}. Then
\be \label{m0}
m_{B}(\bar B, F^{*}(g_{Eucl})) = 0,
\ee
where $m_{B}$ is the Bartnik mass defined by \eqref{bm2}.
\end{theorem} 

\begin{remark}
{\rm The result would be immediate if it were known that the Bartnik mass is continuous (or even lower semi-continuous) 
in the smooth topology on the space of metrics on $\bar B$ of non-negative scalar curvature, since $(\bar B, F^*(g_{Eucl}))$ can be smoothly 
approximated, up to isometry, by domains in $\bR^3$ which have zero Bartnik mass. It is proved in \cite{J2}, \cite{JL} 
that the ADM mass is lower semi-continuous in the pointed $C^{2}$ and $C^{0}$ topologies. However, this does not 
directly imply lower semi-continuity of the Bartnik mass: the main difficulty is that it is not known that ``close'' Bartnik 
data necessarily have ``close''  competitors for near-minimal mass extensions. Note that the recent work of McCormick \cite{Mc3} 
on the continuity of the Bartnik mass does not apply here --- the regions we consider will not generally satisfy the required 
convexity condition in \cite[Theorem 5.1]{Mc3}.
}
\end{remark}

Before proving Theorem \ref{thm_m0}, we first show how it is used to prove Theorem \ref{thm_conj_II}:

\smallskip
  
\noindent {\bf Proof of Theorem \ref{thm_conj_II}:} Consider a pair $(\bar B, F^{*}(g_{Eucl}))$, where $F \in \cF$. Let 
$F_0:S^2 \to \bR^3$ be $F|_{S^2}$, an immersion (but not embedding) of $S^2$ into $\bR^3$. Let $\gamma_0 = F_0^* g_{Eucl}$ and $H_0:S^2 \to \bR$ 
be the induced metric and mean curvature. 

Now, suppose the Bartnik mass of $(\bar B, F^{*}(g_{Eucl}))$ is realized by an extension $(M, g) \in \accentset{\circ} \cP(M)$, so that 
\eqref{bcont_weak} holds with $\Omega = \bar B$. By Theorem \ref{thm_m0}, the ADM mass of $(M,g)$ vanishes. 

Glue $(\bar B, F^{*}(g_{Eucl}))$ and $(M,g)$ along their boundaries so as to satisfy \eqref{bcont_weak}, producing a Riemannian 
manifold $(N,h)$, without boundary, that is asymptotically flat and smooth with non-negative scalar curvature away from the 
gluing hypersurface. The ADM mass of $(N,h)$ also vanishes. By the rigidity case of the positive mass theorem ``with corners'' in 
dimension three \cite{Mi1}, \cite{ST}, $(N,h)$ is isometric to $(\bR^3, g_{Eucl})$ and $H_0 = H_{\partial M}$. In particular, 
there is an isometric \emph{embedding} $G$ of $(\bar B, F^{*}(g_{Eucl}))$ into $(\bR^3, g_{Eucl})$. If we set $G_0 = 
G|_{S^2}$, then $G_0^* g_{Eucl} = F_0^* g_{Eucl} = \gamma_0$ and the mean curvature of the embedding $G_0$ is 
$H_0$. Thus, we have two immersions $F_0$ and $G_0$ of $S^2$ into $\bR^3$ both realizing the same induced metric 
and the same mean curvature. The contradiction will arise because $F_0$ is not an embedding but $G_0$ is, and we will 
show below that $F_0$ and $G_0$ are in fact congruent.

  A pair of immersions $F_{1}$, $F_{2}$ of a surface into $\bR^{3}$ with the same induced metric and mean curvature 
is called a Bonnet pair, and it is well-known that there are no non-trivial Bonnet pairs of spherical topology. We recall the 
simple proof. Let $A_{i}$ be the $2^{\rm nd}$ fundamental form of $F_{i}$. Since $H_{1} = H_{2}$ and $\g_{1} = \g_{2} \equiv \g$, 
the Gauss--Codazzi equations give 
$$\d_{\g}(A_{1} - A_{2}) = 0,$$
where $\d_{\g}$ is the divergence. Also $\tr_{\g}(A_{1} - A_{2}) = 0$. Thus $A_{1} - A_{2}$ is a holomorphic quadratic differential 
on $S^{2}$. Since the only such is $0$, one has $A_{1} = A_{2}$. It then follows from the fundamental theorem for surfaces in 
$\bR^{3}$ (rigidity) that the immersions $F_{1}$ and $F_{2}$ are congruent. 

{\endproof}

\begin{corollary}
The space of compact regions $(\O, g_{\O})$ of non-negative scalar curvature that admit a mass-minimizing extension 
in $\accentset{\circ} \cP(M)$ is not closed in the smooth topology. 
\end{corollary} 

  This Corollary will make it hard to prove the existence of mass-minimizing extensions by studying limits of 
mass-minimizing sequences in general.

\begin{remark}
\label{rem_strict_pos}
{\rm  Recall a result of Huisken--Ilmanen \cite{HI} on the rigidity of the Bartnik mass: if $m_{B}(\O) = 0$, then 
$\O$ is locally flat. (Although note their proof applies to the ``outward-minimizing'' version of the Bartnik mass.) The proof of Theorem \ref{thm_conj_II} above implies the converse of this result is false, for domains $\O$ 
for which Conjecture II holds. To see this, consider a locally flat domain $(\O, g_{\O})=(\bar B, F^*(g_{Eucl}))$, where $F: \bar B \to 
\bR^{3}$ is a smooth immersion which is not an embedding. If $m_{B}(\O) = 0$ {\em and} $m_{B}(\O)$ is realized by a 
minimum-mass extension (i.e.~Conjecture II holds at $\O$), then the proof of Theorem \ref{thm_conj_II} above gives a 
contradiction. Hence, either $m_{B}(\O) > 0$ or Conjecture II fails at $\O$ (or both). 
}
\end{remark}

\noindent {\bf Proof of Theorem \ref{thm_m0}:} Let $F: \bar B \to \bR^{3}$ be the immersion in $\cF$. We detail the proof 
when the self-intersection set $Z \subset S^{2} = \partial \bar B$ consists of two distinct points $z$, $z'$ with $F(z) = F(z')$ 
but with $F$ injective on $\bar B \setminus Z$. The proof in the general case of a finite number of double points 
is a straightforward modification of this case.

  Let $F_0$ denote the restriction $F|_{S^2}$, an immersion (but not embedding) of $S^2$ into $\bR^3$, and let 
$$(\g_{0}, H_{0}) = (F_{0}^*g_{Eucl}, H_{F_{0}})$$ 
be the induced metric and mean curvature of $F_{0}$, defined on $S^2$.

Let $\e > 0$. We will prove that $(\bar B, F^*(g_{Eucl}))$ admits an admissible extension in $\cP(M)$ (for $\e$ sufficiently small) whose ADM mass is 
$\leq C\e$ for a constant $C$ depending only on $F$. For $\e > 0$ small enough, the extension will be in $\accentset{\circ} \cP(M)$, i.e.
it will contain no immersed minimal surfaces surrounding the boundary.
The proof is rather long, consisting of five steps.

\vspace{2mm}
\paragraph{\emph{Step 1: Modification of $g_{Eucl}$ to introduce positive scalar curvature.}}

Let $\Omega \subset \bR^3$ be the image of $F$, a compact set. Fix a neighborhood $U$ of $\Omega$, and let 
$R_0 = \frac{1}{\e}$ be chosen so that the interior of the ball $B(R_0)$ of radius $R_0$ contains the closure of 
$U$ (decreasing $\e$ if necessary).

To begin, we smoothly deform the Euclidean metric $g_{Eucl}$ on $\bR^3$ to produce a new Riemannian metric 
$\hat g$ with the following properties: $\hat g = g_{Eucl}$ on $U$; $\hat g$ has non-negative scalar $\hat s$ that is 
strictly positive somewhere, zero on $\bR^3 \setminus B(R_0)$, and$\int_{\bR^3} \hat s \hat{dv} \leq \epsilon$; 
$\hat g$ is asymptotically flat with ADM mass $\leq \epsilon$. To be definite, we construct $\hat g$ by applying a 
conformal factor $w^4$ to $g_{Eucl}$, where $w$ is superharmonic, harmonic outside $B(R_0)$,  identically 1 
on $U$, and approaches a constant less than 1 at infinity. Additionally, we can choose $w$ so that
\be
\label{w_bound}
1-\epsilon \leq w  \leq 1+ \epsilon, \ \ \ \  1-\epsilon \leq \frac{1}{w}  \leq 1+ \epsilon, \ \ \ \ |\nabla w| \leq \epsilon,
\ee
where $|\nabla w|$ is taken with respect to $g_{Eucl}$. In particular, there exists a closed ball $K$, contained in 
$B(R_0) \setminus U$, such that 
\be \label{a0}
\hat s \geq \alpha_0 > 0 \ \ {\rm on} \ \ K,
\ee
for some constant $\alpha_0 = \a_{0}(\e)$.

\vspace{2mm}
\paragraph{\emph{Step 2: Construction of family of metrics $\bar g_t$ to obtain correct boundary metric.}}
In this step we will perturb the immersion $F_0$ to an embedding. This will of course alter the boundary data 
$(\g_{0}, H_{0})$, so we will also perturb the metric $\hat g$ so as to restore the original boundary metric $\g_{0}$. This 
change will introduce a small amount of negative scalar curvature and will possibly violate \eqref{bcont_weak}; these issues 
will be addressed in Step 3.

\begin{figure}[ht]
\begin{center}
\includegraphics[scale=0.55]{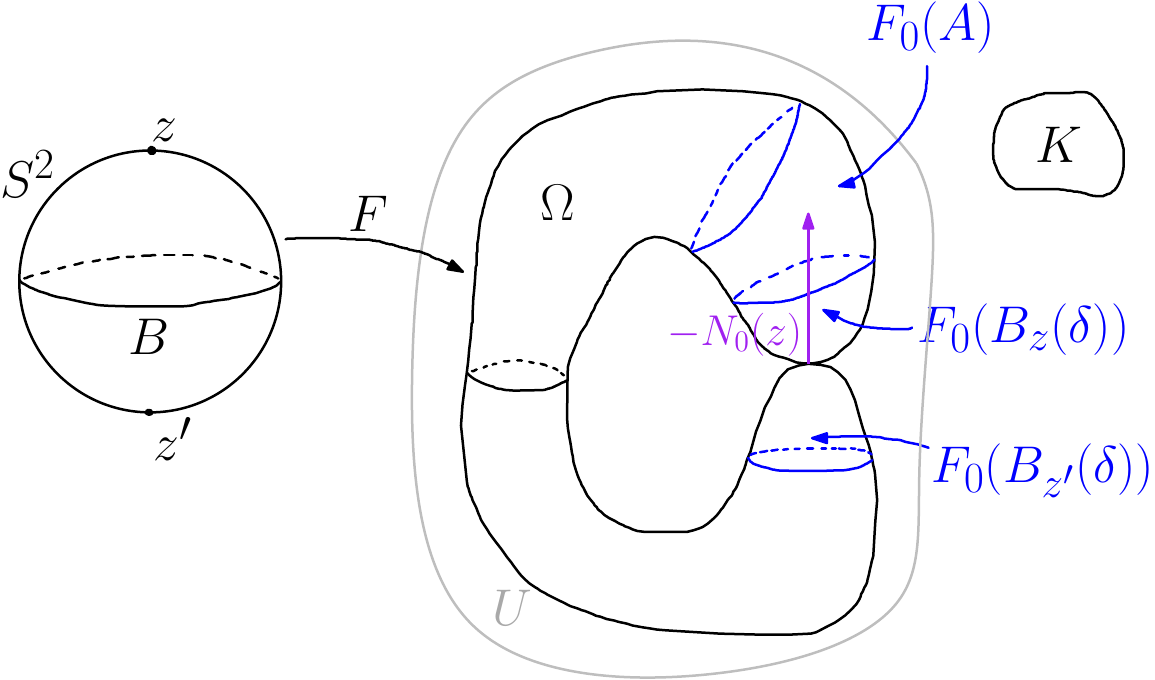}
\end{center}
\caption{\small An illustration indicating the initial setup in the proof of Theorem \ref{thm_m0}.
\label{fig_setup} }
\end{figure}

  To begin, let $N_0$ be the unit outward normal vector field along $F_0$, viewed as a function on $S^2$ (taken with 
respect to $\hat g$, or equivalently, with respect to $g_{Eucl}$). Fix a number $\delta>0$ sufficiently small so that 
$B_z(2\delta) \cap B_{z'}(2\delta) = \emptyset$; here $B_{z}(r) \subset S^2$ is the open geodesic $r$-ball about $z$, 
with respect to the induced metric $\g_{0}$. Let $q: S^{2} \to \bR^{\geq 0}$ be a smooth, non-negative bump function that 
equals $1$ on $B_{z}(\d)$ and is zero outside $B_{z}(2\d)$. Let $A \subset S^2$ be the open annular region:
\be
\label{eqn_A}
A = \interior\left(B_{z}(2\d) \setminus B_{z}(\d) \right).
\ee
An example illustrating this setup is sketched in Figure \ref{fig_setup}.

For $t \geq 0$, define a smooth family of maps $F_t: S^2 \to \bR^3$ by
\be \label{Ft}
F_{t}(x) = F(x) -  tq(x)N_0(z),
\ee
where $N_0(z)$ is treated as a constant vector field on $\bR^{3}$.  For $0 < t < t_{0}$ sufficiently small, $F_t$ is an embedding 
and $F_t(S^2)$ is contained inside $U$. The mapping $t \mapsto F_t(S^2)$ gives a local flow of surfaces in which the set 
$F_t(B_z(\delta))$ is translated in the $-N_0(z)$ direction at speed 1, and $F_t(S^2 \setminus B_z(2\delta))$ does not move. 
Thus, the only change to the geometry occurs in $F_t(A)$.

  For $0< t < t_0$, the smooth, embedded 2-sphere $F_t(S^2)$ bounds a smooth, compact region $\Omega_t$ 
in $\bR^3$ that is diffeomorphic to a closed 3-ball. Let $M_{t} = \bR^{3}\setminus \interior(\O_{t})$, a smooth manifold with 
compact boundary $\dm_{t}$. Note that $F_t$ is a diffeomorphism of $S^2$ onto $\dm_t$. 

\begin{lemma}
\label{lemma_deform_metric}
There exists a domain $V \subset \subset U \subset \bR^3$ and a smooth family 
of Riemannian metrics $\{\bar g_t\}_{0 \leq t < t_0}$ on $\bR^3$ such that:
$$\bar g_0 = \hat g \ {\rm on } \ \bR^3; \ \ \ \bar g_t = \hat g \ \ {\rm outside \ of}  \ V, \  {\rm for} \  t \in [0,t_0),$$ 
and, on $S^{2}$,  
$$F_t^* \bar g_t = F_0^* g_{Eucl} = \g_{0}, \  {\rm for} \  t \in [0,t_0).$$ 
Moreover, $V$ satisfies
\be
\label{eqn_V}
V \cap F_{t}(S^{2}) = F_{t}(\interior(B_{z}({\tfrac{9}{4}}\d) \setminus B_{z}({\tfrac{3}{4}}\d))).
\ee
\end{lemma}

\begin{figure}[ht]
\begin{center}
\includegraphics[scale=0.55]{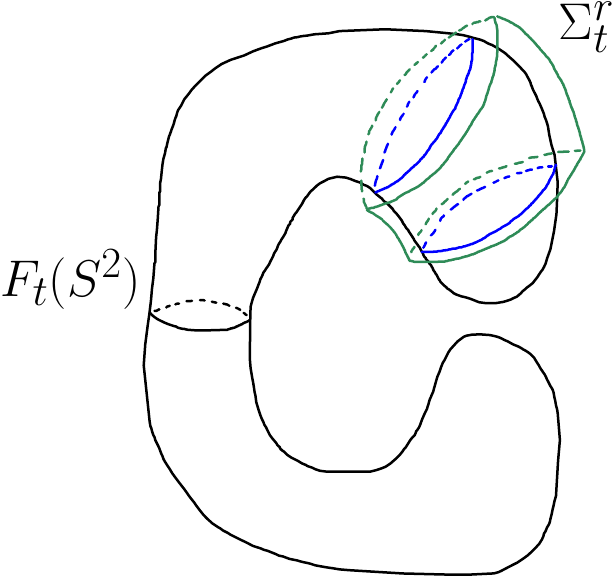}$\qquad\qquad$\includegraphics[scale=0.55]{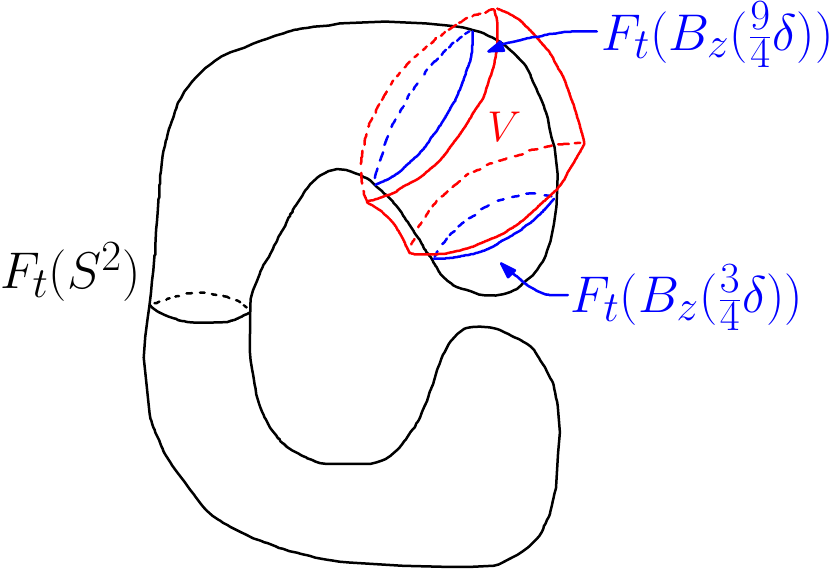}
\end{center}
\caption{\small This illustration shows part of the setup of the proof of Lemma \ref{lemma_deform_metric}.
\label{fig_sigma} }
\end{figure}

\noindent {\bf Proof:} The following construction takes place within $U$, so we regard $g_{Eucl}$ and $\hat g$ as equal in the remainder 
of this proof. Fix $t \in (0,t_0)$, and let $\g_{t}$ be the metric on the embedded surface $F_{t}(S^{2}) \subset \bR^{3}$ induced by 
$g_{Eucl}$.  Let $r(x) = \dist_{g_{Eucl}}(x, \dm_{t})$ be the Euclidean distance of $x \in \bR^3$ to $\dm_{t} = F_{t}(S^{2})$. 
Recall that $\g_{t} = \g_{0}$ outside  $F_{t}(B_z(2\delta))$ (where here we are identifying $\gamma_0 = F_0^*g_{Eucl}$ with 
the induced metric on $F_0(S^2)$). Let $\S_{t}^{r}$ be the image of $F_{t}(A)$ under the time $r$ Euclidean exponential map 
normal to $F_{t}(A)$ (the $r$-equidistant surface to $F_{t}(A)$); here $r \in (-r_{0}, r_{0})$ and $r_{0}$ is chosen small enough 
so that $\S_{t}^{r} \subset M_{t} \cap U$ for $r \geq 0$ and $\S_t^r$ is smooth for all $r \in (-r_0,r_0)$.  Note that $r_{0}$ may 
be chosen independent of $t$; it depends only on $\d$ and the surface $F_{0}(S^{2})$. (See the left side of Figure \ref{fig_sigma}).

 Let $O_t$ be the union of the surfaces $\S_{t}^{r}$ for $|r| < r_{0}$, an open set. In $O_t$, by the Gauss Lemma for the normal 
exponential map, 
$$g_{Eucl} = dr^{2} + \g_{t}^{r},$$
where $\g_{t}^{r}$ is the metric induced on $\S_{t}^{r}$ by $g_{Eucl}$. Note $\g_{t}^{0} = \g_{t}$ on $\Sigma^0_t$. 
Define a new metric on $O_t$ by 
$$\bar g_{t} = dr^{2} + \bar \g_{t}^{r},$$
where $\bar \g_{t}^{r}$ is a smooth metric on $\S_{t}^{r}$, varying smoothly in $r$, such that $\bar \g_{t}^{0} = 
\left.\g_{0}\right|_{\Sigma^0_t}$ and $\bar \g_{t}^{r} = \g_{t}^{r}$, for $|r| \geq r_{0}/2$. (If $t=0$, define 
$\bar g_0 = g_{Eucl}$ on $O_0$.) Note that $\bar \g_t^0$, a metric on $F_t(A)$, extends smoothly to $\gamma_t$ on 
$F_t(S^2)$, since $\gamma_t=\gamma_0$ outside $F_t(A)$.

To complete the construction, let $V$ be an open set contained in $U$ and containing $\displaystyle \bigcup_{0 < t < t_0} O_t$ 
and satisfying \eqref{eqn_V}. (See the right side of Figure \ref{fig_sigma}). 

Extend $\bar g_t$ smoothly to $V$ (and smoothly in $t \in [0,t_0)$ as well, which can be 
arranged in the above construction) so that $\bar g_t$ induces the metric $\gamma_0$ on $F_t(S^2)$ and $\bar g_t$ 
agrees with $\hat g$ near $\partial V$. Then $\bar g_t$ extends to a smooth family of metrics on $\bR^3$, with 
$\bar g_t = \hat g$ outside of $V$. 

{\endproof}

  Note that $\bar g_t = \hat g = g_{Eucl}$ on a neighborhood of $F_0(z)=F_0(z')$. Since the family $\bar g_t$ 
is smooth in $t$ and $\bar g_0 = \hat g$ has non-negative scalar curvature, the scalar curvature $\bar s_t$ of 
$\bar g_t$ is non-negative outside $V$ and converges uniformly to zero inside $V$ as $t \to 0$.  We also note, for later reference, that
$\bar V$ does not contain $F_0(z)$.

  To summarize at this point, we have the asymptotically flat Riemannian manifold $(M_t, \bar g_t)$, for each 
$t \in (0,t_0)$, for which the induced metric on the boundary (when pulled back to $S^2$ via $F_t$) equals 
the original boundary metric $\g_{0}$. The mean curvature $\bar H_t$ of $\partial M_t$ (viewed as a function 
on $S^2$) converges uniformly to the original $H_0$ as $t \to 0$. Moreover, 
\be
\label{barH_H}
\bar H_t = H_0 \qquad \text {on } S^2 \setminus A'
\ee
for all $t \in (0,t_0)$, where $A' = F_{t}^{-1}(V\cap F_{t}(S^{2})) = \interior(B_{z}({\tfrac{9}{4}}\d) \setminus B_{z}({\tfrac{3}{4}}\d))$ 
is a slightly enlarged annulus.

   The space $(M_t, \bar g_t)$ is almost, but not quite, an admissible extension of $(\bar B, F^* (g_{Eucl}))$, for two reasons: first, the mean 
curvatures of the boundaries do not agree (although they are close, for $t$ small) and more generally do not necessarily satisfy 
(\ref{bcont_weak}) (i.e., $\bar H_t \leq H_0$ may fail to hold at some points); second, the scalar curvature $\bar s_{t}$ is not non-negative (although it is nearly 
so, for $t$ small.). We address these two problems in the next step.

\vspace{2mm}
\paragraph{\emph{Step 3: Conformal deformation to correct scalar curvature and boundary mean curvature.}}
Next, we perform a conformal deformation on $(M_{t}, \bar g_t)$.
For each $t \in (0,t_0)$, consider the linear elliptic problem 
\be
\label{eqn_v_t}
\begin{cases}
\bar L_t v_t = 0 & \text{ in } M_t\\
v_t =1 & \text{ on } \partial M_t\\
v_t \to 1 & \text{ at infinity},
\end{cases}
\ee
where $\bar L_t := \bar\D_t - {\tfrac{1}{8}}\bar s_t$ and $\bar\D_t$ is the Laplacian for $\bar g_t$. Assuming for the moment a smooth, positive solution $v_t$ exists, define the 
conformal metric
$$\w g_t = v_t^{4}\bar g_t \qquad \text{on } M_t.$$
Then the induced metric on the boundary stays the same: 
$$F_t^* \w g_t =  F_t^* \bar g_t= F_0^* g_{Eucl} = \g_{0},$$
by Lemma \ref{lemma_deform_metric} and the boundary condition on $v_t$. Also, the mean curvature $\w H_t$ of $\partial M_t$ with respect to $\w g_t$ is given by 
\be \label{wHH1}
\w H_t = \bar H_t + 4\bar N_t(v_t),
\ee
where $\bar N_t$ is the unit  boundary normal on $(M_{t}, \bar g_{t})$ (viewed as a function on $\partial M_t$, 
pointing into $M_t$).  In this step we'll prove that a solution $v_t>0$ to \eqref{eqn_v_t} exists and that, moreover,
\be \label{Hequal}
\w H_t < H_0,
\ee
for $t$ sufficiently small.

  One small difficulty is that, in contrast to the setting of Proposition \ref{prop_H}, since $\bar s_{t}$ may be negative at some points, 
$-\bar L_{t}$ may not automatically be a positive operator (with Dirichlet boundary conditions), so that the equation 
\eqref{eqn_v_t} may not a priori always be uniquely solvable. Similarly, the associated Green's function and 
Poisson kernel may not be uniquely defined, or have appropriate signs. On the other hand, these properties 
are relatively simple to prove. Note first that $\bar L_{t}: C^{k+2,\alpha}_{\d}(M_t) \to C^{k,\alpha}_{\d+2}(M_t)$ 
(where we recall the weighted H\"older space notation from Section 2) is formally $L^2(M_t, \bar g_t)$-self-adjoint 
with respect to zero Dirichlet boundary conditions on $\dm_{t}$. In the statement below, let $M_0 \subset \bR^3$ be the 
closure of the complement of $\Omega$. Note that $M_0 \setminus F_0(z)$ is a smooth manifold with (non-compact) boundary.

Before stating Lemma \ref{pos}, we need a precise definition of a family of functions on $M_t$ converging to a function on $M_0$. 
First, if $\{f_t\}_{0 < t < t_0}$ are in $C^{k,\alpha}_{\d+2}(M_t)$ and $f_0 \in C^{k,\alpha}_{\d+2}(M_0)$, we define $f_t \to f_0$ in $C^{k,\alpha}_{\d+2}(M_0)$ as $t \to 0$ if $f_{t}$ and $f_0$ admit $C^{k,\a}_{\d+2}$ extensions (say $\tilde f_t$ and $\tilde f_0$, respectively)
to a domain $M_* \subset \mathbb{R}^3$ strictly containing 
$M_{0}$ and all $M_t$ (for $t$ sufficiently small) in its interior and $\tilde f_{t} \to \tilde f_{0}$ in $C^{k,\a}_{\d+2}(M_*)$. Second, we define $f_t \to f_0$ in $C^{k,\alpha}_{loc, \d+2}(M_0 \setminus F_0(z))$ if (i) given any compact set $K \subset M_0 \setminus F_0(z)$, there exists a smooth family of embeddings $\Phi_t: K \to M_t$, $\Phi_0$ being the identity, such that $f_t \circ \Phi_t$ converges to $f$ in $C^{k,\alpha}(K)$, and (ii) $f_t |_{\mathbb{R}^3 \setminus U} \to f |_{\mathbb{R}^3 \setminus U}$ in $C^{k,\alpha}_{\delta+2}(\mathbb{R}^3 \setminus U)$.

\begin{lemma} \label{pos}
For $t_0> 0$ sufficiently small, $-\bar L_{t}$ is a positive operator for $0 < t < t_0$, with respect to Dirichlet boundary conditions. 
Hence, given $f_{t} \in C^{k,\alpha}_{\d+2}(M_t)$, there is a unique solution $\omega_{t}\in C^{k+2,\alpha}_{\d}(M_t)$ to 
$\bar L_{t}\omega_{t} = f_{t}$ with $\omega_{t} = 0$ on $\dm_{t}$. Additionally, given $f_{0} \in C^{k,\alpha}_{\d+2}(M_0)$, there is a unique solution $\omega_{0}\in C^{k+2,\alpha}_{\d}(M_0)$ to 
$\bar L_{0}\omega_{0} = f_{0}$ with $\omega_{0} = 0$ on $\dm_{0}$. Moreover, if $f_t \to f_0$ in $C^{k,\alpha}_{\d+2}(M_0)$ as $t \to 0$, then the solutions $\omega_{t}$ converge to $\omega_0$ in 
$C^{k+2,\alpha}_{loc, \d}(M_0 \setminus F_{0}(z))$ as $t \to 0$. Finally, the Green's function $\bar G_{t}(x,y)$ and Poisson kernel $\bar P_{t}(x,y)$ for $\bar L_{t}$ exist and satisfy
$$\bar G_{t}(x, y) \leq 0, \ \ \bar P_{t}(x, y) \geq 0,$$
with strict inequality for $y$ in the interior of $M_{t}$. 
\end{lemma}

\noindent {\bf Proof:} Note that $M_0$ is not a manifold with boundary (as $\partial M_0 = F_0(S^2)$ is not embedded), but it does
satisfy the Poincar\'e ``exterior cone condition" (cf.~\cite[p.29,203--205]{GT}, and is therefore a regular domain for the Dirichlet 
problem for the operator
$$\bar L_{0} = \hat L := \hat \D - {\tfrac{1}{8}}\hat s,$$
where we recall that $\hat s \geq 0$. Clearly, $-\bar L_{0}$ is a positive operator with respect to Dirichlet boundary conditions, 
$\omega = 0$ on $\dm_{0}$ and $\omega \to 0$ at infinity. In particular, for the bottom of the $L^{2}$ spectrum 
one has 
$$\l_{0} = \inf \frac{\int_{M_{0}}-f\bar L_{0}f \overline{dv}_0}{\int_{M_{0}} f^{2}\overline{dv}_0} > 0,$$
where the $\inf$ is taken over nonzero smooth functions $f$ of compact support in $M_{0}$. 
It is standard that there exists a positive Green's function $\bar G_0$ for $\bar L_0$ on $M_0$, and moreover the Dirichlet problem $\bar L_0 \omega_0 = f_0$ is uniquely solvable.

As $t \searrow 0$, the boundaries $\partial M_{t}$ converge to $\partial M_0$, smoothly away from $F_{0}(z)$. Similarly, the operators 
$\bar L_{t}$ converge smoothly to the operator $\bar L_{0}$ away from $F_{0}(z)$ (and are equal outside a compact set). Moreover, the 
lowest eigenvalue $\l_{0}^{t}$ of $-\bar L_{t}$ varies continuously with $t$ as $t \to 0$, (cf.~\cite{Ar} for instance for a much more 
general result than this special situation), and hence $\l_{0}^{t} > 0$, for $t$ sufficiently small. For convenience, we include a direct proof 
that $\l_{0}^{t} > 0$ in the next paragraph.

For $t>0$ let
$$\l_{0}^{t} = \inf \frac{\int_{M_{t}}-f\bar L_{t}f \overline{dv}_t}{\int_{M_{t}} f^{2}\overline{dv}_t} $$
(where the infimum is taken over nonzero smooth functions $f$ of compact support in $M_{t}$), be the bottom eigenvalue of $-\bar L_{t}$. We claim
that $\l_0^t>0$ for $t$ sufficiently small. Since the minimum value of $\bar s_t - \hat s$ on $M_t$ converges to 0 as $t \to 0$,
\begin{align*}
\l_{0}^{t} &= \inf \frac{\int_{M_{t}}\left( |\nabla f|^2_{\bar g_t} + \frac{1}{8} \hat s f^2 \right)\overline{dv}_t + \frac{1}{8} \int_{M_t}\left(\bar s_t - \hat s\right)f^2\overline{dv}_t}{\int_{M_{t}} f^{2}\overline{dv}_t}\\
&= \inf \frac{\int_{M_{t}}\left( |\nabla f|^2_{\bar g_t} + \frac{1}{8} \hat s f^2 \right)\overline{dv}_t }{\int_{M_{t}} f^{2}\overline{dv}_t} - O(t)\\
&\geq c_1 \inf \frac{\int_{M_{t}}\left( |\nabla f|^2_{\hat g} + c_2 \hat s f^2 \right)\hat{dv} }{\int_{M_{t}} f^{2}\hat{dv}} - O(t)
\end{align*}
for some constants $c_1, c_2>0$,  This is because $\{\bar g_t\}_{0 < t < t_0}$ and $\hat g$ are uniformly equivalent by uniform positive constants as $t \to 0$. Now take a manifold with boundary $M_* \subset \mathbb{R}^3$ that contains all $M_t$, and consider $(M_*, \hat g)$.  We have
\begin{align*}
\l_{0}^{t} &\geq c_1 \inf \frac{\int_{M_{*}}\left( |\nabla f|^2_{\hat g} + c_2 \hat s f^2 \right)\hat{dv} }{\int_{M_{*}} f^{2}\hat{dv}} - O(t),
\end{align*}
where this infimum is taken over all nonzero smooth functions $f$ of compact support in $M_*$. It is clear this infimum is strictly positive, so that $\l_0^t>0$ for $t$ sufficiently small.

Thus $-\bar L_{t}$ is a positive operator for $t$ sufficiently small; it is then standard, cf.~\cite{Kr} for instance, that the 
Green's function $\bar G_{t}$ exists and is strictly negative in the interior of $M_{t}$ and hence the Poisson kernel is strictly 
positive (since $\bar P_{t}(x,y) = -N_x\bar G_{t}(x,y)$ and $\bar G_{t}(x,y) = 0$ for $x \in \dm_{t}$). The existence and 
uniqueness, along with the local $C^{k,\alpha}$ convergence  of $\omega_t$ away from $F_{0}(z)$ and the weighted $C^{k+2,\alpha}_{\delta}(\mathbb{R}^3 \setminus U)$ convergence in the end, then follow from standard elliptic estimates. 

{\endproof}

It follows from Lemma \ref{pos} that \eqref{eqn_v_t} has a unique, smooth solution. 

\begin{remark}
\label{remark_P_degen}
{\rm The main technical problem that arises in the discussion to follow is that for $x_{t} \to x = F_{0}(z)$, i.e.~$x_{t}$ 
converging to the singular point, 
$$\bar P_{t}(x_{t}, \cdot) \to 0$$
(uniformly) on $M_{t}$. Thus the Poisson kernel $\bar P_{t}(x_{t}, \cdot)$ degenerates at $F_{0}(z)$. Closely related 
to this is the fact that the Martin boundary of $M_{0}$ equals the Euclidean boundary away from the point 
$F_{0}(z)$ but at the cusp point $F_{0}(z)$ is much larger; there is a minimal positive harmonic function supported at 
$F_{0}(z)$ for each angle of approach to the singular point $F_{0}(z)$. This is discussed in Example 3 of \cite{Ma}. 
}
\end{remark}

  Returning to the analysis of \eqref{eqn_v_t}, it follows as in the discussion concerning \eqref{bulk} that  
\be
\label{eqn_v_t_Green}
v_t(x) = 1 + {\tfrac{1}{8}}\int_{M_t}\bar G_t(x,y)\bar s_t(y)dy,
\ee
for $t \in (0,t_0)$, where $dy$ denotes the volume form $\overline{dv}_t(y)$ of $\bar g_t$. 
As in \eqref{nv}, this gives for $x \in \partial M_t$,
\be \label{nv1}
\bar N_t(v_t)(x) = -{\tfrac{1}{8}}\int_{M_t}\bar P_t(x,y) \bar s_t(y)dy.
\ee

 Let 
$$C_{t}(\d) = F_{t}(B_{z}(\d) \cup B_{z'}(\d)) \ \ \ \ {\rm and} \ \ \ \ D_{t}(\d) = F_{t}(S^{2} \setminus (B_{z}(\d) \cup B_{z'}(\d))).$$
The geometry of $\dm_{t}$ is controlled in $D_{t}(\d)$ but degenerates in $C_{t}(\d)$ as $t \to 0$.
Of course, 
$$\dm_{t} = C_{t}(\d) \cup D_{t}(\d).$$

\begin{lemma}
\label{lemma_Nw0_negative}
If $\epsilon>0$ (from just before Step 1) is sufficiently small, then for $t_0>0$ sufficiently small, the solution $v_t$ to \eqref{eqn_v_t} is positive and satisfies
\be \label{a}
\bar N_t(v_t) < 0, \ \ \ \text{on } \partial M_t, \text{ for } 0 < t < t_0.
\ee
Moreover, there exists $b>0$, independent of $t \in (0, t_0)$, such that
\be \label{b}
\bar N_t(v_t)(x) < -b,  \qquad \text{for } x \in D_{t}({\tfrac{1}{2}}\d),  \ \ 0< t < t_0.
\ee
\end{lemma} 

\noindent {\bf Proof:} By Lemma \ref{pos}, as $t \searrow 0$, $v_t$ converges in $C^{k+2,\alpha}_{loc, \d}(M_0 \setminus F_{0}(z))$ to the limit solution $v_0$ to
\be \label{w0}
\hat L v_0:= \hat \Delta v_0 - \frac{1}{8} \hat s v_0=0,
\ee
on $(M_0, \hat g)$ with boundary conditions $v_{0} = 1$ on $\dm_0$ and $v_{0} \to 1$ at infinity. By the maximum 
principle, since $\hat s \geq 0$ (and is not identically zero), one has the following facts:
\begin{align*} \label{v0}
0< v_0 &\leq 1 \text{ on } M_0,\\
N_0(v_0) &< 0 \text { on }  F_0(S^2 \setminus \{z,z'\}).
\end{align*}
Here, the unit normal $N_0$ is viewed as a (well-defined) vector field on  $F_0(S^2 \setminus \{z,z'\}) \subset \partial M_0$. 
In particular, $N_0(v_0)<-2b$ on $D_0({\tfrac{1}{2}}\d)$ for some constant $b>0$. Since $v_t \to v_0$ locally in $C^1$ away from 
$F_{0}(z)$ (by Lemma \ref{pos}), and $D_t({\tfrac{1}{2}}\d)$ converges smoothly to $D_0({\tfrac{1}{2}}\d)$ as $t \to 0$,
we have
$$\bar N_t(v_t) < -b \text{ on } D_t({\tfrac{1}{2}}\d),$$
for $t$ sufficiently small, which proves \eqref{b}. 

\begin{figure}[ht]
\begin{center}
\includegraphics[scale=0.55]{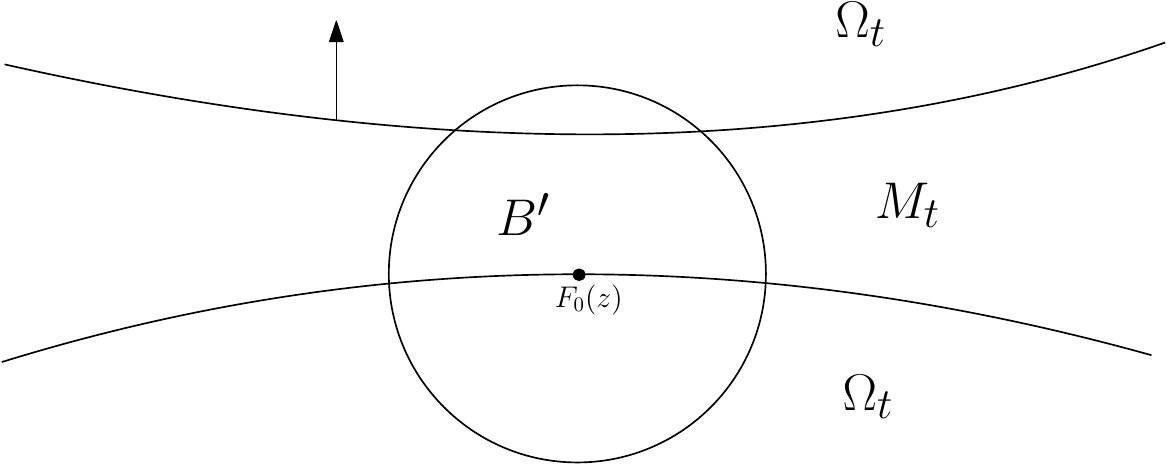}
\end{center}
\caption{\small A sketch of the main sets involved in the proof of \eqref{vtbd}. $B'$ is a fixed ball about $F_0(z)$. $\Omega_t$ is a connected set, but only its ``caps'' are show here, the top one of which translates upward as $t$ increases.
\label{fig_M_t} }
\end{figure}

  Next we claim that $v_{t}$ is uniformly controlled even in a neighborhood of $F_{0}(z)$; specifically we establish \eqref{vtbd} below. Let $B'$ be a Euclidean ball centered at $F_0(z)$ of sufficiently small radius so that $B'$ is disjoint from $\bar V$. Then in the region $M_{t} \cap B'$, the metric $\bar g_{t}$ is the flat Euclidean metric, so that $\bar L_{t}$ 
is the Euclidean Laplacian $\D$; see Figure \ref{fig_M_t}. By the maximum principle applied to the harmonic function $1-v_t$ on the domain $M_t \cap B'$, we have
$$\sup_{x \in M_t \cap B'} |1 - v_{t}(x)| = \sup_{x \in \partial(M_t \cap B')} |1 - v_{t}(x)|.$$
Since $v_t = 1$ on $\partial M_t$, this can be rewritten as
\be\label{sup_v_t}
\sup_{x \in M_t \cap B'} |1 - v_{t}(x)| = \sup_{x \in M_t \cap \partial B'} |1 - v_{t}(x)|.
\ee
By Lemma 4.6, $v_t \to v_0$ in $C^{k,\alpha}_{loc, \d+2}(M_0 \setminus F_0(z))$, which implies (since  $\partial B'$ is a fixed positive distance away from $F_0(z)$), that there exists a compact set $K$ in $M_0$ containing a neighborhood of $M_0 \cap \partial B'$ in $M_0$ and a smooth family of embeddings $\Phi_t: K \to M_t$ such that  $v_t \circ \Phi_t$ converges to $v_0$ in $C^{k,\alpha}(K)$, and that $\Phi_t(K)$ contains a neighborhood of $M_t \cap \partial B'$ (for $t$ sufficiently small). In particular, this implies that $\sup_{x \in M_t \cap \partial B'} |v_0 \circ \Phi_t^{-1} (x)-v_t (x)| \to 0$ as $t \to 0$. From \eqref{sup_v_t} we have
\begin{align*}
\sup_{x \in M_t \cap B'} |1 - v_{t}(x)|  &\leq \sup_{x \in M_t \cap \partial B'} |1 - v_0 \circ \Phi_t^{-1} (x)| +  \sup_{x \in M_t \cap \partial B'} |v_0 \circ \Phi_t^{-1} (x) - v_{t}(x)|,
\end{align*}
and we just argued the second term to the right is $O(t)$. To address the first term on the right, since the metric $\hat g$ is uniformly $\e$-close to the 
Euclidean metric, one has $|1-v_{0}| = O(\e)$ on $M_{t}$.  Thus $\sup_{x \in M_t \cap \partial B'} |1 - v_0 \circ \Phi_t^{-1} (x)|$ is $O(\e)+O(t)$, so
\be \label{vtbd}
\sup_{x \in M_t \cap B'} |1 - v_{t}(x)| \leq O(\epsilon) + O(t),
\ee
as claimed.

Hence, if $\epsilon$ is sufficiently small,  $v_t$ is positive (and more generally $O(\epsilon)$) near $F_0(z)$ for $t$ sufficiently small, proving the claim. In particular, by this together with the appropriate convergence of $v_t$ to $v_0$ as in Lemma \ref{pos}, we have $v_{t} > 0$ on $M_{t}$, for all $t > 0$ sufficiently small.

  The estimate \eqref{a} for $x$ in $C_t(\frac{1}{2}\d)$ for $t$ is somewhat more subtle. Considering \eqref{nv1}, note that 
$\bar P_{t}(x,y) > 0$ and $\bar s_{t} \geq \a_{0} > 0$ in $K$ from \eqref{a0}, while $\bar s_{t}$ slightly negative of order $t$ in 
$V$, for $V$ as in Lemma \ref{lemma_deform_metric}. However, by Remark \ref{remark_P_degen}, the Poisson kernel $\bar P_{t}(x, y)$ 
degenerates at the singular point $x = F_{0}(z)$ as $t \to 0$: $\bar P_{t}(F_{t}(z), y) \to 0$ as $t \to 0$, uniformly in $y$. Thus the 
relative behavior of $\bar P_{t}$ in these two regions is not immediately clear. 

  We will use the following boundary Harnack estimate to obtain uniform control, as $t \to 0$, on the relative behavior of 
$\bar P_t(x,y)$ for $y$ near to and away from the boundary $\dm_{t}$. 

\begin{sublemma}
\label{lemma_harnack}
There exists a constant $C>0$ such that 
\be \label{har}
\sup_{y \in V \cap M_t} \bar P_t(x,y) \leq C \inf_{y' \in K} \bar P_t(x,y'),
\ee
for $x \in C_t(\frac{1}{2}\d)$ and $t \in (0,t_0)$, where $C$ is independent of such $x$ and $t$.
\end{sublemma}
The proof appears later. 

  Now it follows from \eqref{nv1} and the fact that any negative scalar curvature of $\bar g_t$ lies within $V$ together with 
the lower bound $\bar s_{t} \geq \a_{0}$ on $K$ that for $x \in F_t(S^2)$,
$$\bar N_t(v_t)(x) \leq -{\tfrac{1}{8}} \int_{V\cap {M_t}} \bar P_t(x,y)\bar s_t(y)dy -{\tfrac{1}{8}} \int_K \bar P_t(x,y)\alpha_0 dy,$$ 
since $\bar P_t(x,y) \geq 0$.
It follows then from \eqref{har} that for $x \in C_t(\frac{1}{2}\d)$,
\begin{align*}
\bar N_t(v_t)(x) &\leq {\tfrac{1}{8}} \left(\sup_{y \in V \cap M_t} \bar P_t(x,y)\right)\left( \sup_{y \in V \cap M_t} |\bar s_t(y)|\right) \vol_{\bar g_t}(V \cap M_t)  
-{\tfrac{1}{8}} \alpha_0 \left(\inf_{y' \in K} \bar P_t(x,y')\right) \vol_{\hat g}(K)\\
&\leq {\tfrac{1}{8}} \left(\inf_{y' \in K} \bar P_t(x,y')\right)\left[C\left( \sup_{y \in V \cap M_t} |\bar s_t(y)|\right) \vol_{\bar g_t}(V \cap M_t)  
- \alpha_0  \vol_{\hat g}(K)\right].
\end{align*}
Since $\sup_{y \in V \cap M_t} |\bar s_t(y)|$ converges to $0$ (because $\hat s = 0$ on $V$) and $\vol_{\bar g_t}(V \cap M_t)$ is bounded as $t\to 0$, the above is negative for $t$ 
sufficiently small, independent of $x$. This completes the proof of Lemma \ref{lemma_Nw0_negative}.

{\endproof}

Now, we explain why \eqref{Hequal} holds. Recall from \eqref{barH_H} that $\bar H_t = H_0$ on $C_t(\frac{3}{4}\delta)$. Thus, 
\eqref{a} and \eqref{wHH1} show $\tilde H_t < H_0$ on $C_t(\frac{3}{4}\delta)$. Also, $\bar H_t$ converges uniformly to $H_0$, 
so \eqref{wHH1}  and \eqref{b} show (shrinking $t_0$ if necessary) that $\tilde H_t < H_0$ on $D_t(\frac{1}{2}\delta)$. This proves 
\eqref{Hequal}. 

To conclude this step, we note that $\w g_t=v_t^4 \bar g_t$ is asymptotically flat: since $\bar s_t$ vanishes outside a compact set, $v_t$ is 
$\bar g_t$-harmonic outside a compact set. Since $v_t \to 1$ at infinity, it is well-known (and not hard to show) that $v_t^4 \bar g_t$ is 
asymptotically flat. Moreover, $\w g_t$ has zero scalar curvature since $\bar L_t v_t=0$. Thus, $\w g_t$ is an admissible extension of 
$(\bar B, F^*(g_{Eucl}))$ in $\cP(M)$, for $t$ sufficiently small.

\vspace{2mm}

\paragraph{\emph{Step 4: Control of ADM mass of $\w g_t$.}}

By the conformal deformation formula \eqref{confmass}, and the fact that $\hat g = \bar g_0 = \bar g_t$ outside a compact set,
\begin{align}
m_{ADM}(\w g_t) &= m_{ADM}(\bar g_t) - \frac{1}{2\pi} \lim_{r \to \infty} \int_{S_r} \bar N_t(v_t) \overline{dA}_t \nonumber\\
&= m_{ADM}(\hat g) - \frac{1}{2\pi} \int_{S_{R_1}} \bar N_t(v_t) \overline{dA}_t  
\leq \epsilon - \frac{1}{2\pi} \int_{S_{R_1}} \bar N_0(v_t) \overline{dA}_0 \nonumber\\
&\leq \epsilon + \frac{1}{2\pi} \int_{S_{R_1}} | \bar \nabla_0 (v_t)| \overline{dA}_0, \label{eqn_mass_est}
\end{align}
by the divergence theorem (since $v_t$ is $\bar g_t$-harmonic outside $U$ for all $t$) and since $\bar g_t = \bar g_0$ on 
$S_{R_1}$. Here $R_1 \geq R_0+1$, where $R_{0}$ is the value chosen in Step 1, i.e.~$S_{R_{0}}$ encloses $U \supset \dm_{t}$. 
Increasing $R_1$ if necessary, we arrange that
\be
\label{S_r}
\frac{|S_{R_1}|_{\bar g_0}}{4\pi (R_1)^2} \leq 2,
\ee
by asymptotic flatness.

By Lemma \ref{pos}, the convergence of $v_t$ to $v_0$ is sufficient to guarantee that, by \eqref{eqn_mass_est}, 
\be
\label{m_tilde}
m_{ADM}(\w g_t)  \leq \epsilon + \frac{2}{2\pi} \int_{S_{R_1}} | \bar \nabla_0 (v_0)| \overline{dA}_0
\ee
for $t$ sufficiently small. Using the Green's function $\bar G_0$ to represent $v_0$, we have, as in 
\eqref{eqn_v_t_Green}:
$$v_0(x) = 1 + {\tfrac{1}{8}}\int_{M_0}\bar G_0(x,y)\bar s_0(y)dy = 1 + {\tfrac{1}{8}}\int_{B(R_0)}\bar G_0(x,y)\bar s_0(y)dy,$$
since $\bar s_0$ vanishes outside $B(R_0)$.
By the standard decay of the Green's function, there exists constant $C_1$ depending only on the initial immersion $F$ such that
$$|\bar \nabla_0 \bar G_0(x,y)|_{\bar g_0} \leq \frac{C_1}{|x|^2}$$
for $|x| \geq R_1$ and $y \in B(R_0)$. From the $\e$-bound on the $L^1$ norm of the scalar curvature of $\bar g_0 = \hat g$ 
from Step 1, this gives 
$$|\bar \nabla_0 v_0(x)|_{\bar g_0} \leq \frac{C_1\e}{8|x|^2},$$
for $|x| \geq R_1$. Combining this with \eqref{m_tilde} and using \eqref{S_r} implies that 
\begin{align*}
m_{ADM}(\w g_t) &\leq \epsilon + \frac{2}{2\pi} \frac{C_1 \epsilon}{8(R_1)^2} |S_{R_1}|_{\bar g_0}\\
&\leq (1+C_1)\epsilon.
\end{align*}

\vspace{2mm}
\paragraph{\emph{Step 5: Absence of Horizons}} In this final step, we argue that if $\e>0$ was chosen small enough to 
begin with in Step 1, then $(M_t, \w g_t)$ will not contain any immersed minimal surfaces that surround $\partial M_t$, 
for $t$ sufficiently small.

Recall from Step 1 that $R_0=\frac{1}{\e}>0$ was chosen so that $B(R_0)$ contains $U$, and $\hat g$ was constructed to be 
conformally flat with with zero scalar curvature outside $B(R_0)$.  In particular, $\bar g_t$ and $\w g_t$ also have zero scalar 
curvature and are conformally flat outside $B(R_0)$. It follows that $\w g_t = u_t^4 g_{Eucl}$ in $\bR^3 \setminus B(R_0)$, 
where $u_t>0$ is $g_{Eucl}$-harmonic for each $t$ (specifically, $u_t = wv_t$).

In particular, using Lemma \ref{pos} and the decay of $w$, one has the Euclidean estimates 
\be \label{Eest}
|u_t(x)-1+a| \leq \frac{c_0}{|x|} \ \ {\rm and} \ \ |\nabla u_t (x)|_{g_{Eucl}} \leq \frac{c_1}{|x|^2},
\ee
for $x \in \bR^3 \setminus B(R_0)$, where $c_0$ and $c_1$ depend only on $F$, and $1-a$ is the constant that $w$ 
(and hence $u_t$) approaches at infinity. By \eqref{w_bound}, $0<a \leq \epsilon$.
 
  Suppose $\S$ is a compact, immersed minimal surface in $(M_t, \w g_t)$ that surrounds $\partial M_t$. The mean curvatures $\w H_{t}$ and 
$H_{Eucl}$ of $\S$ with respect to $\w g_t$ and $g_{Eucl}$ are related by
$$0 = \w H_t = (1+f_1) H_{Eucl} + f_2$$
for smooth functions $f_1$ and $f_2$ on $\S$, where $|f_1|$ and $|f_2|$ are bounded above by the $C^1$ norm of 
$\w g_t-g_{Eucl}$. (This can be seen from the first variation of area formula, for instance). Fix a Euclidean ball $B''$ centered at $F_0(z)$ that does not contain $\Sigma$. In particular, by \eqref{w_bound},  the fact that $\bar g_t$ converges smoothly to $\hat g$ as $t \to 0$, and that $v_t$ converges smoothly to $v_0$ outside $B''$, we have, for $t$ sufficiently small, 
\be
\label{H_Eucl1}
|H_{Eucl}|  \leq C_2 \e \qquad \text{ on } \Sigma \setminus B'',
\ee
for a constant $C_2$ depending only on $F$ and $B''$.

Let $r = |x|$ be the Euclidean distance function from $F_0(z)$, and let $R_{1} = \max_{\S}r$, 
achieved at a point $p_{0} \in \S \setminus B''$. By a standard comparison of mean curvature of $\S$,
\be 
\label{H_comparison}
H_{Eucl}(p_{0}) \geq H_{Eucl}(S_{R_{1}}) = \frac{2}{R_{1}}.
\ee
Thus, $\frac{1}{R_1} \leq \frac{C_2 \epsilon}{2}$. 
If $\epsilon>0$ is sufficiently small, then
$R_1 > R_0$. In particular, $p_0 \in \bR^3 \setminus B(R_0)$, so that $\w g_t = u_t^4 g_{Eucl}$ in a neighborhood of $p_0$.
Thus
$$0 = \w H_t = u_t^{-2} H_{Eucl} + 4 u_t^{-3}N_{Eucl}(u_t),$$
at $p_0$. Combining this with \eqref{Eest},
\be
\label{H_Eucl2}
|H_{Eucl}(p_{0})| \leq 4 \frac{|\nabla u_t|}{u_t}(p_{0}) \leq \frac{4 c_1}{(R_1)^2} \frac{1}{1-\epsilon-\frac{c_0}{R_1}} 
\ee
Estimates \eqref{H_comparison} and \eqref{H_Eucl2} give a contradiction if $\e$ is sufficiently small, since $\e$ controls $\frac{1}{R_1}$.

\vspace{2mm}

Thus, if $\e$ is chosen to be sufficiently small in Step 1, then $\w g_t$ is an admissible extension of $(\bar B, F^* (g_{Eucl}))$ in $\accentset{\circ} \cP(M)$ 
for $t$ sufficiently small. The proof of Theorem \ref{thm_m0} is now complete, except for the proof of Sub-Lemma \ref{lemma_harnack}, 
to which we now return.  

{\endproof}

\noindent {\bf Proof of Sub-Lemma \ref{lemma_harnack}.} 
We will use the boundary Harnack principle for the elliptic operator $\bar L_t$, cf.~\cite[Theorem 1.1]{BB},
 for instance. Recall that the open set $V$, 
from Lemma \ref{lemma_deform_metric}, satisfies $F_0(z) \not\in \bar V$, $V \cap F_t(S^2) = F_t(A') = F_t( \interior(B_{z}({\tfrac{9}{4}}\d) \setminus B_{z}({\tfrac{3}{4}}\d)))$ 
for each $t$, and that $K$ is a set disjoint from $U$ on which $\hat s \geq \alpha_0 > 0$. 

\begin{figure}[ht]
\begin{center}
\includegraphics[scale=0.55]{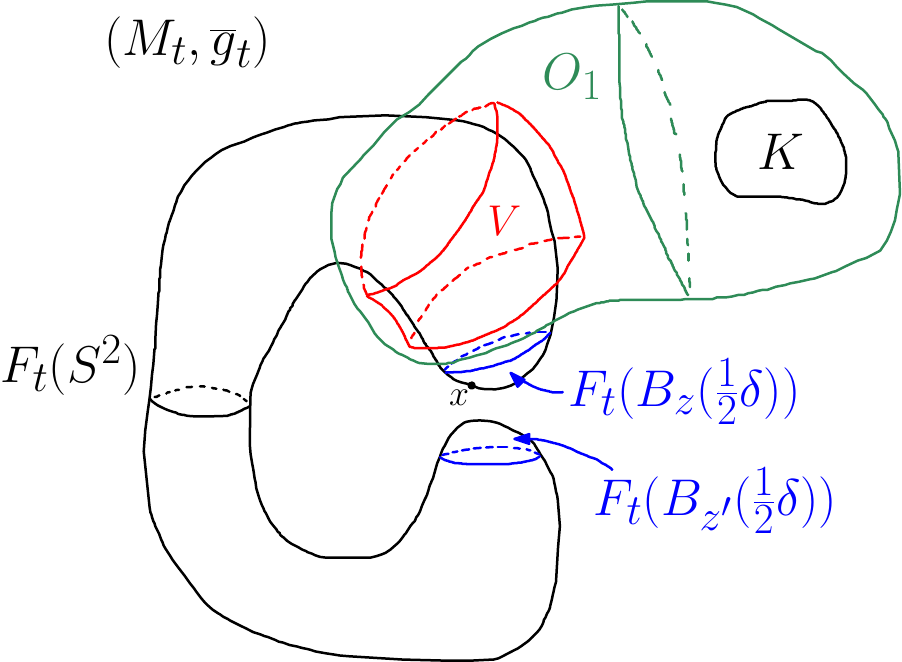}
\end{center}
\caption{\small The main sets used in the proof of Sub-Lemma \ref{lemma_harnack} are shown above. $O_2$ is not pictured but can be 
viewed as a slight enlarging of $O_1$.
\label{fig_harnack} }
\end{figure}

Let $O_2 \supset O_1 \supset (V \cup K)$ be connected, bounded open sets in $\bR^3$ chosen so that 
$\overline O_2 \cap C_t(\frac{5}{8}\delta) = \emptyset$ for all $t \in (0,t_0)$
and that $\overline O_1 \subset O_2$. In particular, $\overline O_2$ does not include $F_0(z)$. See Figure \ref{fig_harnack}.

Then for any $x \in C_t(\frac{1}{2}\delta)$, the function 
$y\mapsto \bar P_t(x,y)$ is smooth and bounded on $O_2 \cap M_t$ and vanishes on $\partial M_t \cap O_2$.

 Let $\nu_t$ by the unique solution to $\bar L_t \nu_t=0$ in $M_t$ with boundary conditions  
$0$ on $\dm_t$ and 1 at infinity. By Lemma \ref{pos}, $\nu_t$ converges smoothly, away from $F_{0}(z)$, as $t \to 0$, 
to a function $\nu_0$ on $M_0$ satisfying $\hat L \nu_0=0$, cf.~also the discussion around \eqref{vtbd}. By the maximum principle, $0 < \nu_0 < 1$ 
in the interior of $M_0$. The convergence in the region $O_2$ (which again is a finite distance from $F_0(z)$) is sufficient to guarantee that, shrinking $t_0$ if necessary,
\begin{itemize}
\item $0< \nu_t \leq 2 $ in $O_2 \cap \interior(M_t)$,
\item $\nu_t(y) \geq \beta$ for $y \in K$, where $\beta>0$ is some constant independent of $t$.
\end{itemize}

  By Lemma \ref{pos} again, for $x \in C_t(\frac{1}{2}\delta)$ and $y \in O_2$, the Poisson kernel $\bar P_{t}(x,y)$ 
satisfies $\bar P_{t}(x, \cdot) = 0$ on $O_2 \cap \dm_{t}$ and $\bar P_{t}(x, \cdot) > 0$ in $O_2 \cap \interior(M_{t})$.   
Since $\bar L_t$ is elliptic and $\nu_t(y)$ and $y \mapsto \bar P_t(x,y)$ are both $\bar L_t$-harmonic on $O_2 \cap M_t$, the 
boundary Harnack principle (cf.~\cite[Theroem 1.1]{BB}) implies that: for all $y,y' \in O_1 \cap \interior (M_t)$ and all 
$x \in C_t({\tfrac{1}{2}}\delta)$, one has 
\be
\label{bhp}
\frac{\bar P_t(x,y)}{\nu_t(y)} \leq c_t\frac{\bar P_t(x,y')}{\nu_t(y')}
\ee
for some constant $c_t>0$ (depending on $t$), but independent of $x$. However, since $\bar L_t$ and $(M_t, \bar g_t)$ 
converge smoothly in the region $O_2$ (away from $F_{0}(z)$) as $t\to 0$, we may take the constant $c_t$ independent 
of $t \in (0,t_0)$; call it $C_0$.

Thus, using \eqref{bhp} and the relations on $\nu_{t}$ above, 
\begin{align*}
\sup_{y \in V \cap M_t} \bar P_t(x,y) &= \sup_{y \in V \cap \interior(M_t)} \bar P_t(x,y)
\leq \sup_{y \in O_1 \cap \interior(M_t)} \bar P_t(x,y) \leq \sup_{y \in O_1 \cap \interior(M_t)} \frac{2 \bar P_t(x,y)}{\nu_t(y)}\\
&\leq 2C_0\inf_{y' \in O_1 \cap \interior(M_t)} \frac{\bar P_t(x,y')}{\nu_t(y')} 
\leq 2C_0\inf_{y' \in K} \frac{\bar P_t(x,y')}{\nu_t(y')}
\leq \frac{2C_0}{\beta}\inf_{y' \in K} \bar P_t(x,y').
\end{align*}
This proves the result with $C= \frac{2C_0}{\beta}$. 

{\endproof} 

The proof of Theorem \ref{thm_m0} is now complete.

\begin{remark}
\label{rmk_bcs}
{\rm In the proof of Theorem \ref{thm_m0}, we constructed admissible extensions of $(\bar B, F^*(g_{Eucl}))$ that obeyed the boundary 
conditions \eqref{bcont_weak}. However, it may be possible with some further work to achieve equality of the mean curvatures, 
i.e.~\eqref{bcont}, in the construction by following an argument similar to the proof of Proposition \ref{prop_H}. Specifically, 
one may replace \eqref{eqn_v_t} with $\bar L_t v_t = f_t$, where the functions $f_t \geq 0$ are chosen to be supported near 
$\dm_t$ and so that the normal derivatives $\bar N_t(v_t)$ satisfy \eqref{wHH1} with $\bar H_t = H_0$.  Note that  $f_t$ may blow up near the singular point $F_0(z)$, which complicates the analysis; we do not pursue this further here.
}
\end{remark}

  To conclude this section, we note that it is not difficult to see that the proof of Theorem \ref{thm_conj_II} generalizes to a larger class of 
immersions $\cF$ at the boundary of the space of embeddings than the particular class $\cF$ used in Theorem \ref{thm_conj_II}. We will 
not pursue this in any further detail here. Instead, we make the following more general:

\begin{conjecture}
\label{conj_false}
Conjecture II is false for any {\it locally flat} 3-ball. That is, if $F$ is any smooth immersion of a 3-ball $\bar B$ in 
$\bR^{3}$ that is not an embedding, then $(\bar B, F^{*}(g_{Eucl}))$ admits no
admissible extension realizing its Bartnik mass. 
\end{conjecture}

\section{Remarks on Conjecture III}
\setcounter{equation}{0}

In this section, we discuss several aspects of Conjecture III, related to the analysis in the previous section on Conjecture II. 

   To begin, (as noted briefly in the Introduction), it is proved in \cite{A1}, \cite{AK} that the moduli space $\cE^{m,\a}$ of $C^{m,\a}$ 
AF static vacuum solutions $(g, u)$, $u > 0$, on $M = \bR^{3}\setminus B$ is a smooth Banach manifold. The moduli space 
$\cE^{m,\a}$ is the space of all such static vacuum metrics $(g, u)$ which are $C^{m,\a}$ smooth up 
to $\dm$, modulo the action of the $C^{m+1,\a}$ diffeomorphisms ${\rm Diff}_{1}^{m+1,\a}(M)$ of $M$ 
equal to the identity on $\dm$ and asymptotic to the identity at infinity. Moreover, the map to Bartnik boundary data 
\be \label{PiB1}
\Pi_{B}: \cE^{m, \a} \to Met^{m,\a}(S^{2})\times C^{m-1,\a}(S^{2}) := \cB,
\ee
$$\Pi_{B}(g,u) = (\g, H),$$
is a smooth Fredholm map, of Fredholm index 0, i.e.~$\dim \Ker D\Pi_{B} = \dim \Coker D\Pi_{B}$, at any $(g, u)$. 

  The ADM mass of $(M, g)$ is given by a simple Komar integral 
\be \label{mf}
m_{ADM}(g) = \frac{1}{4\pi} \int_{\dm}N(u)dA_{\g},
\ee
and clearly the mass
\be \label{m}
m_{ADM}: \cE^{m,\a} \to \bR,
\ee
is a smooth function on $\cE^{m,\a}$. 

  Conjecture III is the statement that the map $\Pi_{B}$ is a bijection when $H > 0$. As in \eqref{Pi+}, let $\cE_{+}^{m, \a}$ 
be the open Banach submanifold of static vacuum metrics with $H > 0$ at $\dm$. The map 
\be \label{piB}
\Pi_{B}: \cE_{+}^{m, \a} \to Met^{m,\a}(S^{2})\times C_{+}^{m-1,\a}(S^{2}), 
\ee
$$\Pi_{B}(g) = (\g, H),$$
is clearly also a smooth Fredholm map, of Fredholm index 0. 
  
  The question of whether $\Pi_{B}$ in \eqref{piB} is a bijection is a PDE issue (global existence and 
uniqueness for an elliptic boundary value problem) which is now disconnected from the extension issue in Conjecture I. 
The interior behavior in $B = \bR^{3}\setminus M$ no longer plays any role (besides assigning boundary data). 
In particular, the mass function $m_{ADM}$ in \eqref{m} may well have negative values on $\cE_{+}^{m,\a}$. Put another way, 
it is not at all clear (at least in general) how to restrict the boundary data $(\g, H)$ to the smaller space $\cB_{+}$ of such data 
which have non-negative scalar curvature in-fillings in order to obtain meaningful information about the restricted map 
$\Pi_{B}|_{\cD_{+}}$, where $\cD_{+} = \Pi_{B}^{-1}(\cB_{+})$. 

\medskip 

  It is proved in \cite{AK} that $\Pi_{B}$ in \eqref{piB} is not a homeomorphism. In fact, $\Pi_{B}$ is not proper, and 
if the inverse map is defined, it is not continuous. The reasons for this are more or less the same as the behavior 
discussed in Theorem \ref{thm_conj_II}, namely the passage from embedded spheres to immersed spheres, and it is 
worth discussing this in more detail. 

  Let $\Imm^{m+1,\a} := \Imm^{m+1,\a}(S^{2}, \bR^{3})$ be the space of $C^{m+1,\a}$ immersions $F: S^{2} \to \bR^{3}$. 
This is a smooth Banach manifold (an open submanifold of the full mapping space $C^{m+1,\a}(S^{2}, \bR^{3})$). Similarly the 
space $\Emb^{m+1,\a} := \Emb^{m+1,\a}(S^{2}, \bR^{3})$ of $C^{m+1,\a}$ of embeddings is an open submanifold of  
$\Imm^{m+1,\a}$.  Of course embeddings $F \in \Emb^{m+1,\a}$ give static vacuum solutions $(M, g_{Eucl}, 1)$ where 
$M$ is the unbounded component of $\bR^{3}\setminus \Image F$; thus 
\be \label{ee}
\Emb^{m+1,\a} \subset \cE^{m+1,\a}
\ee
upon an appropriate identification. Immersions that are not embeddings no longer give such flat static vacuum solutions. 
It is then natural to consider the behavior of the inclusion \eqref{ee} at the (point-set theoretic) boundary of $\Emb^{m+1,\a}$ 
within $\Imm^{m+1,\a}$; denote this space as $\partial \Emb^{m+1,\a}$. 
 
 In the following, we will identify immersions into $\bR^{3}$ that differ by a rigid motion of 
$\bR^{3}$. Rigid motions, i.e.~the isometry group of $\bR^{3}$, act freely on $\Imm^{m+1,\a}$ by 
post-composition. Let $\cI mm^{m+1,\a}$ be the resulting smooth quotient space.

   For $F \in \cI mm^{m+1,\a}$, the induced metric $\g = F^{*}(g_{Eucl})$ is a $C^{m,\a}$ metric on $S^{2}$ while the mean curvature 
$H = H_{F}$ is in $C^{m-1,\a}(S^{2})$. Note that the data $(\g, H)$ are well-defined for $F \in \cI mm^{m+1,\a}$. Thus the map 
$\Pi_{B}$ in \eqref{PiB1} or \eqref{piB}, defined initially on $\Emb^{m+1,\a}$ extends to a smooth map on the larger space 
$\cI mm^{m+1,\a}$; to avoid confusion, we denote this extended map as $\Pi_{B}^{\cI}$.  
   
\begin{lemma} 
The map 
$$\Pi_{B}^{\cI}: \cI mm^{m+1,\a} \to \cB, \ \ \Pi_{B}^{\cI}(F) = (\g, H) = (F^{*}(g_{Eucl}), H_{F}),$$
is a smooth proper embedding of Banach manifolds. 
\end{lemma}

\noindent {\bf Proof:} The map $\Pi_{B}^{\cI}$ is injective by the proof of Theorem \ref{thm_conj_II}, i.e.~the non-existence of 
(non-trivial) Bonnet pairs. The proof that $D\Pi_{B}^{\cI}$ is injective is essentially the same. Thus, the Gauss--Codazzi (constraint) 
equations for the immersion $F$ are $\d_{\g}(A - H\g) = 0$, where $\d_{\g}$ is the divergence. Linearizing gives $\d'_{\g'}(A - H\g) + \d_\g(A' - H'\g - H\g') = 0$, where $\d'_{\g'}$ is the variation of the divergence. If $(\g', H') 
= (0,0)$, this becomes $\d_\g A' = 0$. Since $\tr A' = (\tr A)' + \<A, \g'\> = 0$, it follows as before that $A'$ is a holomorphic quadratic 
differential on $S^{2}$ and hence $A' = 0$. Thus the full Cauchy data $(\g', A')$ of the immersion vanish. It follows by 
(infinitesimal) rigidity of surfaces that $F'$ is an infinitesimal rigid motion, so $F' = 0$ in $T\cI mm^{m+1,\a}$. 

Next we show that $\Pi_{B}^{\cI}$ is proper. Suppose $F_{i}$ satisfy $\Pi_{B}^{\cI}(F_{i}) = (\g_{i}, H_{i}) \to (\g, H)$ in 
$\cB$. Then one has uniform control on $\d A$ and $\tr A$. It is well-known that $(\d, \tr)$ form an elliptic system for 
symmetric bilinear forms on $S^{2}$. Since the system has trivial kernel on $\cI mm^{m+1,\a}$, elliptic regularity gives 
uniform control on $\{A_{i}\}$ in $C^{m-1,\a}$. It is then standard that this gives uniform control on $\{F_{i}\}$ in $\cI mm^{m+1,\a}$. 
Thus a (sub)-sequence of $\{F_{i}\}$ converges to a limit $F$, which proves that $\Pi_{B}^{\cI}$ is proper.

{\endproof}

 Let 
$$\cM = {\rm Im}(\Pi_{B}^{\cI}) \subset \cB,$$
a properly embedded Banach submanifold representing the Bartnik boundary data of immersions 
$F: S^{2} \to \bR^{3}$. Let 
$$\cM^{emb} \subset \cM$$
be the open submanifold of embedded Bartnik boundary data, i.e.~$\cM^{emb} = \Image(\Pi_{B}^{\cI}(\Emb^{m+1,\a}))$. 
Thus 
$$\cM^{emb} \subset \Image(\Pi_{B}).$$
However, it is not at all clear if the full space of immersed boundary data $\cM \subset \Image(\Pi_{B})$. 

  As discussed in \cite{AK}, the map $\Pi_{B}$ is not proper when restricted to $\Emb^{m+1,\a} \subset \cE^{m,\a}$. 
Namely, take any sequence of embeddings $F_{i}$ converging to an immersion $F$ that is not an embedding. The Bartnik boundary data 
$(\g_{i}, H_{i})$ of $F_{i}$ converge to the boundary data $(\g, H)$ of $F$. However, the sequence of static vacuum solutions 
$(M, g_{i}, 1)$ determined by $F_{i}$ does not converge to a limit in $\cE^{m,\a}$ (since $\partial(\Emb^{m+1,\a})$ is not 
contained in $\cE^{m,\a}$). 

\medskip 
 
In analogy to Conjecture \ref{conj_false}, we make: 

\begin{conjecture} 
\label{conj_III_false}
{\rm The Bartnik boundary data $(\g, H)$, with $H > 0$, of any locally flat 3-ball (that is not an embedded ball in $\bR^3$) 
does not have a static vacuum extension, i.e.~such $(\g, H) \in \cM \setminus \cM^{emb}$ are not in the image of $\Pi_{B}$. 
}
\end{conjecture}

  One may similarly conjecture that Conjecture \ref{conj_III_false} above holds more generally for immersions $F: S^{2} \to M$ into 
any static vacuum solution $(M, g, u)$ in place of flat $\bR^{3}$, where $F(S^{2})$ is a surface surrounding $\dm$. 

\medskip 

  We also conjecture there is a second region where Conjecture III breaks down. Recall that the black hole uniqueness 
theorem \cite{I}, \cite{BM}, together with \cite{Mi2}, states that the only static vacuum extension of the 
boundary data $(\g, 0)$ is given by the Schwarzschild metric with $\g = \g_{2m}$ a round metric of radius $2m$. 

\medskip 

\begin{conjecture}
{\rm  For any $\g \in Met^{m,\a}(S^{2})$ of non-constant Gauss curvature, 
there is a neighborhood $\cU_{\g} \subset C_{+}^{m-1,\a}(S^{2})$ with $0 \in \overline \cU_{\g}$, such that 
for $H \in \cU_{\g}$, the boundary data $(\g, H)$ does not bound a static vacuum metric $(M, g, u)$. 
In particular, $\Pi_{B}$ is not surjective near the Schwarzschild metric. 
}
\end{conjecture}

   Partial evidence for this conjecture is given by the main compactness theorem (Theorem 1.2) in \cite{AK}. Namely, 
if $(M, g_{i}, u_{i})$ is a sequence of static vacuum solutions, $u_i>0$, with boundary data $(\g_{i}, H_{i}) 
\to (\g, H)$ in $\cB$, with $H_{i} > 0$ and $H \geq 0$, and if $\dm$ is strictly outer-minimizing in $(M, g_{i}, u_{i})$ for all $i$, 
then a subsequence of $(M, g_{i}, u_{i})$ converges in $\cE^{m,\a}$ to a limit $(M, g, u)$ realizing the data $(\g, H)$. 
Setting $H = 0$, one has a contradiction to the black hole uniqueness theorem if $\g \neq \g_{2m}$ for 
some $m > 0$. Hence, such outer-minimizing solutions cannot exist for $i$ sufficiently large. 

  It remains an open question as to whether such static vacuum extensions exist with $\dm$ not 
outer-minimizing. If such a sequence exists, either the curvature of $g_{i}$ must blow-up near $\dm$ or the distance 
to the cut-locus of the normal exponential map must tend to zero (or both), as $i \to \infty$. 

\begin{remark}
{\rm 
  The compactness result above suggests modifying the Bartnik mass $m_{B}$ to $\tilde m_{B}$ by allowing for only 
outer-minimizing extensions, as suggested by Bray \cite{Br}. Note this rules out the constructions above 
in the proof of Theorem \ref{thm_conj_II} and the discussion above on Conjecture III, which are certainly {\em not} 
outer-minimizing extensions. One may also restrict the map $\Pi_{B}$ in \eqref{piB} to the space $\w \cE^{m,\a}$ of 
static vacuum solutions for which $\dm$ is strictly outer-minimizing. Note that $\w \cE^{m,\a}$ is an open domain in 
$\cE_{+}^{m,\a}$. However, as discussed above, the restricted map $\Pi_{B}$ on $\w \cE^{m,\a}$ is not 
surjective (onto a product neighborhood of the Schwarzschild boundary data). Thus, Conjecture III also 
fails for the modified mass $\w m_{B}$. 

  Observe that boundary data $(\g, H)$ near Schwarzschild data $(\g_{2m}, 0)$ {\em do} have outer-minimizing 
extensions in $\accentset{\circ} \cP^{m,\a}$. The discussion above (together with Theorem \ref{thm1.1}) strongly suggests that Conjecture II also fails for the modified 
mass $\w m_{B}$, i.e.~there exist $(\g, H)$ for which there is no mass-minimizing extension realizing $\w m_{B}$.

}
\end{remark}

  Although Theorem \ref{thm_conj_II} shows that Conjecture II is false in general, (and similarly the discussion above 
indicates that Conjecture III is likely to be false in general) one would still like to find natural geometric conditions on the 
boundary data $(\g, H)$ of the region $\O$ under which these Conjectures could remain valid. In a simpler but related setting, 
the guiding light along these lines is the famous Weyl embedding theorem \cite{N}, \cite{P} that a 2-sphere $S^{2}$ with metric 
$\g$ of positive Gauss curvature embeds isometrically in $\bR^{3}$ as the boundary of a convex body $\O$. In particular, the 
normal exponential map $\exp_{N}$ into the exterior $M = \bR^{3}\setminus \O$ has no cut or focal points. 

  As discussed in \cite{AK}, \cite{A2}  and seen here in Theorem \ref{thm_conj_II}, the presence of nearby cut or focal 
points of $\exp_{N}$ is the primary difficulty in establishing Conjecture III and is of course also basic in establishing 
Conjecture II. Thus, it is natural to ask: 

\vspace{1mm}

\noindent {\bf Question.}  Are there natural geometric conditions on $(\g, H)$ such that any extension of $(\g, H)$ in $\cP^{m,\a}$ or 
any static vacuum extension in $\cE_{+}^{m,\a}$, has a lower bound on the distance to the cut-locus of $\exp_{N}$? 

\vspace{2mm}

  Unfortunately, there is little evidence (if any) to suggest that the conditions $K_{\g} > 0$ and $H > 0$ are sufficient 
for this purpose, i.e.~a simple, direct generalization of the Weyl embedding theorem has little support for its validity. 
On the other hand, it would of course be interesting to find any examples where $K_{\g} > 0$, $H > 0$ with the 
distance to the cut-locus of $\exp_{N}$ arbitrarily small.

\bibliographystyle{plain}

\end{document}